\documentclass{article}
\usepackage{latexsym,amsmath,amssymb}
\usepackage{color}
\usepackage[colorlinks]{hyperref}
\usepackage{graphics}            
\usepackage{graphicx}
\usepackage{mathrsfs}
\usepackage[all]{xy}

\newtheorem{definition}{Definition}[section]

\newtheorem{lemma}[definition]{Lemma}
\newtheorem{proposition}[definition]{Proposition}
\newtheorem{theorem}[definition]{Theorem}

\newtheorem{remark}[definition]{Remark}

\newcommand{\proof}{\textbf{Proof : }}
\date{}
\begin{document}
\title{Extension of local positive definite $\mathbb{Z}_{2}^{n}$-superfunctions 
}
\author{Fatemeh Nikzad Pasikhani,$^{\rm
a}$\footnote{
Email: fatemeh.nikzad@iasbs.ac.ir; f.nikzad.p@gmail.com} \ \
Mohammad Mohammadi$^{\rm a}$\footnote{E-mail: m.mohamady64@gmail.com }\\
Saad Varsaie$^{\rm a}$\footnote{E-mail: varsaie@iasbs.ac.ir}
 \\$^a$Department of Mathematics, Institute for Advanced Studies in \\ Basic Sciences, Zanjan, Iran.
}

 \maketitle \glossary{}{}
\begin{abstract}
	By defining a concept of  boundedness for $ \mathbb{Z}_2^n $-superfunction, we show that each local positive definite $ \mathbb{Z}_2^n $-superfunctions, which is bounded in some sense, on a $ \mathbb{Z}_2^n $-Lie supergroup has a positive definite extension to all of the $\mathbb{Z}_{2}^{n}$-Lie supergroup.
\end{abstract}

\maketitle

\section{Introduction}
The theory of positive definite functions is important in different parts of pure and applied mathematics, including harmonic analysis, functional analysis, spectral theory, representation theory, representation of Lie groups, operator theory. it is also used in fields such as mathematical statistics, approximation theory, optimization, quantum physics and probability models, such as stochastic processes.
The theory of positive definite functions is also related to the representation of the unitary group, for more details, see \cite{Sasvari1944,Berg2008}. 
In the simple case, positive definite functions on  additive group $\mathbb{R}$ are important in the study of spectral theory for Schrodinger equations. See \cite{JorgensenPedersen2016,Stewart1976}.

\medskip

In 1940, Krein \cite{Krein1940} proved that every scalar-valued positive definite function on the interval $(-a,a)$ where $a$ is a positive real number has a positive definite continuous extension to the entire real line.
He used an extension theorem for positive linear functionals to prove that this extension always exists but is not necessarily unique.
A direct proof is provided in \cite{Raikov1940}. 
In 1945, Livshic \cite{Livshic1945} studied the problem of extension a continuous positive definite function on $(-2a,2a) \times \mathbb{R}$ to a positive definite function on the plane.
He used operator methods in his proof.
Akutovich \cite{Akutowicz1959} and Devinatz \cite{Devinatz1959} explained the uniqueness of the extension of a positive definite scalar-valued function on an interval in terms of a self-adjoint operator.
The later author introduced the necessary conditions for existence of extensions of such functions when the domain is a rectangle $(-a_1,a_1) \times (-a_2,a_2) \subset \mathbb{R}^2$. 
In 1963, Rudin \cite{Rudin1963} showed that for $ n \geq 2 $ there exist continuous positive definite functions on the rectangle of $\mathbb{R}^n$ that do not extend to a continuous positive definite on $\mathbb{R}^n$.
Rudin extended Krein's result under additional assumptions and proved the following:
Every radial positive definite function defined on a ball in $\mathbb{R}^n$ around zero extends to a radial positive definite function defined on the whole space $\mathbb{R}^n$.
Therefore, we must apply some symmetry conditions on the function or its domain in order to guarantee the existence of an extension.
On Lie groups, generally, locally defined  positive definite functions do not extend to positive definite functions on the whole group.
By considering additional assumptions, this extension is possible; see for example \cite{Devinatz1959,Jorgensen1989,Jorgensen1990,Jorgensen1991,Rudin1970}.
Friedrich \cite{Friedrich1991} in 1991 investigated the results and questions related to the problem of extending a continuous positive definite function from an open neighborhood of a neutral element of a Lie group to the entire group.
In 2001, the extension of a special class of locally definite positive functions on the Heisenberg group was studied in \cite{Bruzual2001}.
Also the problem of the extension of continuous positive definite  functions on closed subgroups of a nilpotent locally compact group was studied in \cite{Kaniuth2004}.

\medskip

In 2011, Neeb \cite{Nee11} proved that every locally definite positive analytic function on certain simple connected Frechet-Lie groups expands to an analytic function on the whole group.
In 2012, Neeb and Salmasian \cite{NeebSalmasian2013} proved that for a broad class of Frechet-Lie supergroups there is a correspondence between smooth (resp., analytic) positive definite superfunctions and the matrix coefficients of smooth (resp., analytic) unitary representations of the Harish-Chandra pair related to this supergroups.
In 2017, Mohammadi and Salmasian \cite{MMSH} extend the theory of unitary representations to the $\mathbb{Z}_{2}^{n}$-graded case
and show that there is a correspondence between smooth positive definite superfunctions and cyclic unitary representations.

\medskip

In this article we show that
local positive definite uniformly bounded  (c.f. Definition \ref{e6}) $ \mathbb{Z}_2^n $-superfunctions extend on the whole $ \mathbb{Z}_2^n $-Lie supergroup. For this, we need to prove stability theorem (c.f. Theorem \ref{f2}) and to construct a pre-representation corresponding to a positive definite $\mathbb{Z}_{2}^{n}$-superfunction for the $\mathbb{Z}_{2}^{n}$-Harish-Chandra pair of the $\mathbb{Z}_{2}^{n}$-Lie-supergroup. In addition we extend the pre-representation to a smooth unitary representation by using stability theorem. It is worth to mention that we show that stability theorem holds unconditionally. Our proof is different from what is stated in \cite{MMSH}.

\medskip

The outline of this article is as follows. In Section \ref{e14}, some basic concepts and preliminaries to $ J $-adically Hausdorff complete and $\mathbb{Z}_{2}^{n}$- supergroups are stated. In Section \ref{e13}, some concepts are introduced such as $\mathbb{Z}_{2}^{n}$-Hilbert spaces, pre-representati-ons and smooth unitary representations. In Section \ref{e12}, we prove the Stability Theorem on $\mathbb{Z}_{2}^{n}$- Lie supergroups.
In Section \ref{e9},
we introduce the concepts of uniformly boundedness  and local positive definiteness for $\mathbb{Z}_{2}^{n}$-superfunctions. In Section \ref{e10}, we construct a pre-representation corresponding to a positive definite $\mathbb{Z}_{2}^{n}$-superfunction of a $\mathbb{Z}_{2}^{n}$-Harish-Chandra pair. In Section \ref{e11}, it is shown that each local positive definite $\mathbb{Z}_{2}^{n}$-superfunction can be
extended if it is uniformly bounded.

\textbf{Acknowledgements.} We would like to express our deep and sincere gratitude to Professor Hadi Salmasian from the University of Ottawa for proposing this problem and providing invaluable guidance throughout this research.
\section{Preliminaries}\label{e14}
This section includes two subsections. By referring to \cite{covolo2016differential}, In the first subsection, we have given  the concept of $ J $- adically Hausdorff completeness. In the second subsection $\mathbb{Z}_{2}^{n}$- supergroups have been discussed. See \cite{MMSH,mohammadi2021construction} for more details.
\subsection{$ J $-adically Hausdorff complete}
Let $ J $ be a homogeneous ideal of a $\mathbb{Z}_{2}^{n}$-ring $ R $, i.e.
$J=\bigoplus\limits_{m \in \mathbb{Z}_{2}^{n}}(J\cap R^{m})$
and $ M $ be an $ R $-module. The collection of sets $ \{x + J^kM\}_{k=0}^{\infty} $ is called the $ J $-adic topology, i.e. it is a basis for a topology on $ M $, where $ x $ runs over all elements of $ M $. It is clear
from the definition that the $ J $-adic topology is translation-invariant with respect to the additive group
structure of $ M $. Now if $ f : M \rightarrow N $ is an $ R $-module morphism. Then
$ f $ is $ J $-adically continuous.
\begin{definition}
	Let $J$ be an ideal of $R$. $R$ is called $J$-adically Hausdorff if for any two arbitrary points $x$ and $y$ of $R$, there exist $J$-adically open neighborhoods $x+J^n$ and $y+J^n$ of $x$ and $y$ respectively such that for an $ n $, their intersection is empty. 
\end{definition}
\begin{proposition}
	Let $J$ be an ideal of comutative ring $R$. $R$ is $J$-adically Hausdorff if and only if $ \bigcap\limits_{n=1}^{\infty}J^{n}=0 $.
\end{proposition}
\begin{definition}
	Let $J$ be an ideal of $R$. A sequence of elements of $ R $, say $ \{a_i\} $, is a $ J $-adic Cauchy sequence if for each non-negative integer $ n $ there exists $ l $ such that $ a_i - a_k \in J^n $ for all $ i, k \geq l $. $ a_j $
	converges $ J $-adically if there exists $ a \in R $ such that for each non-negative integer $ n $, there exists $ l $ such
	that $ a_j - a \in J^n $ for all $ j \geq l $.
\end{definition}
\begin{definition}
	Let $J$ be an ideal of $R$. $R$ is called $ J $-adically complete if every $ J $-adically Cauchy sequence in $ R $ converges $ J $-adically to a unique limit in $ R $.
\end{definition}
\begin{theorem}
	Let $R$ be a commutative $\mathbb{Z}_{2}^{n}$-superalgebra and $J$ a homogeneous ideal. $R$ is $ J $-adically Hausdorff complete if and only if the natural ring morphism $P: R \rightarrow \varprojlim\limits_{k \in \mathbb{N}}
\frac{R}{J^k}$ is an isomorphism.
\end{theorem}
\begin{definition}
	Let $X$ be a topological space and $\mathcal{O}_{X}$ a sheaf of $\mathbb{Z}_{2}^{n}$-comutative $\mathbb{Z}_{2}^{n}$-ring. A pair $(X,\mathcal{O}_{X})$ is called a $\mathbb{Z}_{2}^{n}$-ringed space. Let $(X,\mathcal{O}_{X})$ and $(Y,\mathcal{O}_{Y})$ be $\mathbb{Z}_{2}^{n}$-ringed spaces. A pair $\psi:=(\overline{\psi},{\psi}^*)$ is a morphism of them if $\overline{\psi}: X \rightarrow Y$ is a continuous map and $\psi^*:\mathcal{O}_{Y} \rightarrow \overline{\psi}_{*}\mathcal{O}_{X}$ is a morphism of weight $0$.
\end{definition}
In what follows, we assume that $ n $-tuples $ \gamma_{0},\gamma_{1},...,\gamma_{2^{n}-1}$
are elements of $\mathbb{Z}_2^n $ in  lexicographical order: $\gamma_{0}\prec \gamma_{1}\prec...\prec \gamma_{2^{n}-1}$, where $\gamma_{0}=0$. Let $\gamma_s=\sum_{j=1}^{n}a_je_j$ and $\gamma_t=\sum_{j=1}^{n}b_je_j$, the lexicographical order  $\gamma_s \prec \gamma_t $ means that $a_{j}=\bar{0}$, $b_{j}=\bar{1}$ for some $1 \leq j \leq n$ and $a_{k}=b_{k}$ for all $k < j$.\\
Let $ \gamma_s=(a_1,a_2,...,a_n), \gamma_t=(b_1,b_2,...,b_n) \in \mathbb{Z}_{2}^{n}$, we define $\mathcal{B}:\mathbb{Z}_{2}^{n} \times \mathbb{Z}_{2}^{n} \rightarrow \{\pm1\}$ as follows 
$$\mathcal{B}(\gamma_s,\gamma_t):=(-1)^{\langle \gamma_s,\gamma_t \rangle}\enspace
$$
where $ \langle \gamma_s,\gamma_t \rangle=a_1b_1+a_2b_2+...+a_nb_n.$
\begin{definition}
	Let $n,p,q_1,q_2,...,q_{2^{n}-1} \in \mathbb{N}$ and	$\mathfrak{q}=(q_{1},q_{2},...,q_{2^{n}-1}) $. 
	Assume $p$ coordinates $x^1,x^2,...,x^p$ of degree $\gamma_{0}=0$, $q_{1}$ coordinates 	$ \xi^{1},\xi^{2},...,\xi^{q_{1}} $ of degree
	$ \gamma_{1} $,
	$ q_{2} $
	coordinates
	$ \xi^{q_{1}+1},\xi^{q_{1}+2},...,\xi^{q_{1}+q_{2}} $
	of degree
	$ \gamma_{2} $,...and 
	$q_{2^{n}-1}$
	coordinates
	$ \xi^{q_{1}+q_{2}+...+q_{2^{n}-2}+1},...,\xi^{q_{1}+q_{2}+...+q_{2^{n}-2}+q_{2^{n}-1}} $
	of degree
	$ \gamma_{2^{n}-1} $. Set  $x=(x^{1},x^{2},...,x^{p})$ of all coordinates of degree zero and $\xi=(\xi^{1},...,\xi^{q})$ of all coordinates of degree nonzero such that $q=\sum\limits_{k\geq 0}q_{k}$. A $\mathbb{Z}_{2}^{n}$-ringed space
	\begin{align*}
	&&\mathbb{R}^{p|\mathfrak{q}}:=(\mathbb{R}^p,C_{\mathbb{R}^{p}}^{\infty}(-)[[\xi^1,...,\xi^q]])&&
	\end{align*}
	is called a $\mathbb{Z}_{2}^{n}$-superdomain of dimension $p|\mathfrak{q}$ where $C_{\mathbb{R}^{p}}^{\infty}$ is sheaf of smooth functions over $\mathbb{R}$. Also for every open set $U \subset \mathbb{R}^{p},$
	\begin{align*}
	&&\mathcal{O}_{\mathbb{R}^{p|\mathfrak{q}}}(U):=C_{\mathbb{R}^{p}}^{\infty}(U)[[\xi^1,...,\xi^q]]&&
	\end{align*}
	it is $\mathbb{Z}_{2}^{n}$-commutative associative unital $\mathbb{Z}_{2}^{n}$-algebra of formal power series with indeterminates or formal variables $\xi^{a}$ of degree $\gamma_a$ in	$\mathbb{Z}_{2}^{n} \backslash \{0\}$ and its members are in the form of a formal power series as follows
	\begin{align*}
	&&f(x,\xi)=\sum\limits_{\lvert\mu\rvert=0}^{\infty} f_{\mu_1...\mu_q} (\xi^1)^{\mu_1}...(\xi^q)^{\mu_q}=\sum\limits_{\lvert\mu\rvert=0}^{\infty} f_\mu(x) \xi^{\mu}&&
	\end{align*}
	where coefficients are in the ring $C_{\mathbb{R}^{p}}^{\infty}(U)$ and multiplication is subject to
	the sign rules determined by 
	$$\xi^i \xi^j=\mathcal{B}(\gamma_i,\gamma_j) \xi^j \xi^i.$$
\end{definition}
\begin{proposition}
	Algebra $\mathcal{O}(V)=C^{\infty}(V)[[\xi^{1},...,\xi^{q}]]$ of $\mathbb{Z}_{2}^{n}$-functions over $V$ is  $ J $-adically Hausdorff complete.
\end{proposition}
\subsection{$\mathbb{Z}_{2}^{n} $-Lie supergroups}
\begin{definition}
	Let $p \in \mathbb{N}$ and $\mathfrak{q}=(q_1,...,q_{2^{n}-1}) \in \mathbb{N}^{2^{n}-1}$. A $\mathbb{Z}_{2}^{n}$-(smooth) supermanifold 
	$\mathcal{M}$
	of dimension 
	$p|\mathfrak{q}$ is a 
	$\mathbb{Z}_{2}^{n}$-local ringed space 
	$(M,\mathcal{A}_M)$
	which is locally isomorphic to 
	$\mathbb{Z}_{2}^{n}$-superdomain
	$ (\mathbb{R}^{p},C_{\mathbb{R}^p}^{\infty}[[\xi^1,...,\xi^q]]) $
	where
	$q=\sum_{k} q_k$
	and
	$\xi^1,...,\xi^q$ are
	$\mathbb{Z}_{2}^{n}$-commutative
	formal variables and
	$q_k$
	is the number of formal variables which  their degrees are equal to $ k $.
	and
	$C_{\mathbb{R}^p}^{\infty}$ is the sheaf of the smooth real functions of Euclidean space
	$\mathbb{R}^p$.
\end{definition}
Let $ \mathbb{R}^{0|0} $ be
a  $\mathbb{Z}_{2}^{n}$-superdomain with a singleton $\{0\}$ as a topological space and the constant sheaf $\mathbb{R}$ on $\{0\}$. It is a terminal object in the category of $\mathbb{Z}_{2}^{n}$-supermanifold. Indeed, for each $\mathbb{Z}_{2}^{n}$-supermanifold $ \mathcal{M}=(M,\mathcal{A}_M)  $, there exists a unique morphism $\mathcal{I}_{M}=(I_{M},I_{M}^*):\mathcal{M} \rightarrow \mathbb{R}^{0|0}$ as follows:
\begin{align*}
\begin{array}{ccc}
I_M:& M \rightarrow \{0\}, \qquad\qquad I_M^*:\mathbb{R} \rightarrow \mathcal{A}_M,\\ & \qquad\qquad\qquad\qquad\qquad k \mapsto k. 1_{\mathcal{A}_M} 
\end{array}
&
\end{align*}
In addition for each $p \in M$, there is a morphism $\mathfrak{j}_p=(j_p,j_p^*):\mathbb{R}^{0|0}\rightarrow \mathcal{M}$ that $j_p:\{0\} \rightarrow M$ and
\begin{flalign*}
\begin{array}{cc}
j_p^*:& \mathcal{A}_M \rightarrow \mathbb{R}\\ &  g \mapsto \tilde{g}(p):=ev_p(g).
\end{array}
&
\end{flalign*}
Therefore, for every $\mathbb{Z}_{2}^{n}$-supermanifold $\mathcal{T}$, the following constant morphism can be defined
\begin{align*}
&&\hat{P_T}:\mathcal{T} \rightarrow \mathbb{R}^{0 |0} \overset{\mathfrak{j}_p}{\rightarrow}\mathcal{M}&&
\end{align*}
It is the composition of $\mathfrak{j}_p$ and the unique morphism $\mathcal{T} \rightarrow \mathbb{R}^{0|0}$. 
Therefore, for each $f \in \mathcal{A}_{M}$ we have
\begin{align*}
&& \hat{P_T}:\mathcal{T} \rightarrow \mathcal{M}, \qquad \hat{p_T}^{*}: f \mapsto ev_{p}(f) . 1_{\mathcal{A}_M}&&
\end{align*}
Therefore the category of $\mathbb{Z}_{2}^{n}$-supermanifold has a terminal object $ \mathbb{R}^{0|0} $.
\begin{definition}
	Let $G$ be a $\mathbb{Z}_{2}^{n}$-real smooth supermanifold. G is called a $\mathbb{Z}_{2}^{n}$-real Lie supergroup if there exist the following morphisms
	\begin{itemize}
		\item[i)] $\mu: G \times G \rightarrow G$ (multiplication morphism)
		\item[ii)] $i:G \rightarrow G$ (inverse morphism)
		\item[iii)] $e:\mathbb{R}^{0|0} \rightarrow G$ (unit morphism) where $ \mathbb{R}^{0|0} $ is a terminal object in  the category of $\mathbb{Z}_{2}^{n}$-supermanifold. 
	\end{itemize}
	such that following diagrams are commutative
	\begin{displaymath} 
	\scalebox{0.8}{\xymatrix{
			& G\times G \times G \ar[rd]^{\mu\times 1_{G}} \ar[ld]_{1_{G} \times \mu} & \\
			G\times G \ar[dr]_\mu &  & G\times G \ar[ld]^{\mu}\\
			& G & }\quad
		\xymatrix{
			& G\times G \ar[rd]^{\mu} & \\
			G\ar[ur]^{<<1_{G},\hat{e}_{_G}>>} \ar[dr]_{<<\hat{e}_{_G},1_{G}>>}  \ar[rr]_{1_{G}} &  & G\\
			& G\times G \ar[ur]_\mu  & }\quad}
	\scalebox{0.8}{\qquad	\xymatrix{
			& G\times G \ar[rd]^{\mu} & \\
			G\ar[ur]^{<<1_{G},i>>} \ar[dr]_{<< i,1_{G}>>}  \ar[rr]_{\hat{e}_{_G}} &  & G\\
			& G\times G \ar[ur]_\mu  & }}
	\end{displaymath}
	where $\hat{e}$ denotes the composition of the identity map $ e $  with the unique map
	$ G \rightarrow \mathbb{R}^{0|0} $.
	Moreover, $ <<\varphi,\psi>> $ denotes the map $(\varphi \times \psi) \circ \mathrm{d}_{G}$ where $\mathrm{d}_{G}: G \rightarrow G\times G$ is the canonical diagonal map.
\end{definition}
\begin{definition}
	A $\mathbb{Z}_{2}^{n}$-superalgebra $\mathfrak{g}=\bigoplus\limits_{a \in \mathbb{Z}_{2}^{n}} \mathfrak{g}_a$ is caled  $\mathbb{Z}_{2}^{n}$-Lie superalgebra if there exists $\mathbb{Z}_{2}^{n}$-superbracket
	$[.,.]: \mathfrak{g} \times \mathfrak{g}\rightarrow \mathfrak{g}$ on $\mathfrak{g}$ such that
	\begin{itemize}
		\item[i)] $[. , .]$ is bilinear and for every $a,b \in \mathbb{Z}_{2}^{n}$, we have  $[\mathfrak{g}_a, \mathfrak{g}_b] \subset \mathfrak{g}_{a+b}$, 
		\item[ii)] for every  $x \in \mathfrak{g}_a$ and $y \in \mathfrak{g}_b$ where $a,b \in \mathbb{Z}_{2}^{n}$
		$$[x,y]=- \mathcal{B}(a,b)[y,x]\enspace (Graded\enspace skew-symmetry),$$
		\item[iii)] for every
		$x \in \mathfrak{g}_a$,
		$y \in \mathfrak{g}_b$
		and
		$z \in \mathfrak{g}_c$
		where
		$ a,b,c \in \mathbb{Z}_{2}^{n}$ we have
		\begin{align*}
		&&[x,[y,z]]=[[x,y],z]+\mathcal{B}(a,b) \mathcal{B}(a,c)[y,[z,x]]\enspace(Graded\enspace Jacoby\enspace identity).&&
		\end{align*}
	\end{itemize}
\end{definition}
If $g$ is a  $\mathbb{Z}_{2}^{n}$-Lie superalgebra, then its  0-degree part ,$\mathfrak{g_{0}}$, is a classical Lie algebra.\\
Let $\mathfrak{g_{\mathbb{C}}}:=\mathfrak{g}\otimes_{\mathbb{R}}\mathbb{C}$. We define $\mathfrak{U}(\mathfrak{g_{\mathbb{C}}})$ the universal enveloping algebra of $ \mathfrak{g_{\mathbb{C}}} $ as follows:
$$\mathfrak{U}(\mathfrak{g_{\mathbb{C}}}):=\dfrac{T(\mathfrak{g}_{\mathbb{C}})}{I}$$
where $ T(\mathfrak{g}_{\mathbb{C}}) $ denotes the tensor algebra of $ \mathfrak{g_{\mathbb{C}}} $ in the category
$\mathbb{Z}_{2}^{n}$-graded complex vector spaces
and $ I $ denotes the two-sided ideal of $ T(\mathfrak{g}_{\mathbb{C}}) $ generated by elements of the form
$$x \otimes y-\mathcal{B}(a,b)y\otimes x-[x,y]$$
where $x,y$ are homogeneous elements in $\mathfrak{g}_{\mathbb{C}}$ of degrees $ a,b $ respectively. See \cite[Sec.2.1]{NK} for further details.\\
Set $ \textbf{u}(a)=|\{1\leq j\leq n\enspace : \enspace a_j=\bar{1}\}| $
for any $a=(a_1,...,a_n) \in \mathbb{Z}_{2}^{n}$ and define 
\begin{align}
\alpha(a)=e^{\frac{\pi i}{2} \textbf{u}(a)}
\end{align}
Suppose $x \mapsto x^*$ is the (unique) conjugate-linear map on $\mathfrak{g_{\mathbb{C}}}$ defined as follows, for each $x \in \mathfrak{g}_{a}$ and $a \in \mathbb{Z}_{2}^{n}$
\begin{align}\label{e4}
x^*:=-\overline{\alpha(a)}x.
\end{align}
As $I$ is a $*$-invariance,
the map $x \mapsto x^*$ extends to a conjugate-linear antiautomorphism of the algebra $\mathfrak{U}(\mathfrak{g_{\mathbb{C}}})$. 
Then $(D_{1}D_{2})^*=D_{2}^*D_{1}^*$ for every $D_1,D_2 \in \mathfrak{U}(\mathfrak{g_{\mathbb{C}}})$.
\begin{definition}
	Let $G_0$ be a Lie group and $\mathfrak{g}=\bigoplus\limits_{a \in \mathbb{Z}_{2}^{n}} \mathfrak{g}_a$ be a $\mathbb{Z}_{2}^{n}$-Lie superalgebra. A pair $(G_0, \mathfrak{g})$ is called a $\mathbb{Z}_{2}^{n}$-Harish-Chandra pair if there exists an action $Ad:G_0 \times \mathfrak{g} \rightarrow \mathfrak{g} $ such that preserve the  $\mathbb{Z}_{2}^{n}$-grading and $ \mathrm{Ad}|_{\mathfrak{g_{0}}}:G_{0} \times \mathfrak{g_{0}} \rightarrow \mathfrak{g_{0}} $ is adjoint action of $G_0$ on $\mathfrak{g_{0}} \cong Lie(G_{0})$.
\end{definition}
\begin{definition}
	Let $ (G,\mathfrak{g}) $
	and $(H,\mathfrak{h})$ be two $\mathbb{Z}_{2}^{n}$-Harish-chandra pairs. A morphism $(\phi,\varphi):(G,\mathfrak{g}) \rightarrow (H,\mathfrak{h})$ is defined as fallows
	\begin{itemize}
		\item [i)] map $\phi: G \rightarrow H$ is a Lie group homomorphism,
		\item [ii)] map $ \varphi: \mathfrak{g}\rightarrow \mathfrak{h} $ is a morphism in category of $\mathbb{Z}_{2}^{n}$-Lie superalgebras such that $ \varphi|_{\mathfrak{g_{0}}}=\mathrm{d}\phi $,
		\item [iii)] $\mathrm{Ad}_{H} \circ (\phi,\varphi)=\varphi \circ \mathrm{Ad}_{G}$ for $\mathrm{Ad}_{H}$ 
	\end{itemize}
\end{definition}
\begin{proposition}\cite[Pro 2.6]{MMSH}
	The category of $\mathbb{Z}_{2}^{n}$-Lie supergroups is equivalent to the category of $\mathbb{Z}_{2}^{n}$-Harish-Chandra pairs.	
\end{proposition}
Let
$\mathcal{G}:=(G_{0},\mathcal{O}_{G_0})$ denote the $\mathbb{Z}_{2}^{n}$-Lie supergroup corresponding to $\mathbb{Z}_{2}^{n}$-Harish-chandra pair $ (G_0, \mathfrak{g}) $. Assume $C^{\infty}(\mathcal{G})$ is a $\mathbb{Z}_{2}^{n}$-superalgebra of smooth function on $ \mathcal{G} $. Then 
$$ C^{\infty}(\mathcal{G})\cong \mathrm{Hom}_{\mathfrak{g}_0}(\mathfrak{U}(\mathfrak{g}_{\mathbb{C}}),C^{\infty}(G_0,\mathbb{C}))  $$
in which $\mathfrak{g_{0}}$ acts on $ C^{\infty}(G_0,\mathbb{C}) $ as the left invariant differential operators.
\section{smooth unitary representations}\label{e13}
In this section, following \cite{MMSH}, we further deal with our other needed preliminaries.
\subsection{$ \mathbb{Z}_{2}^{n} $-Hilbert superspace}
For any $a=(a_1,...,a_n) \in \mathbb{Z}_{2}^{n}$, we recall that
$\alpha(a)=e^{\frac{\pi i}{2} \textbf{u}(a)}$, where $ \textbf{u}(a)=|\{1\leq j\leq n\enspace : \enspace a_j=\bar{1}\}| $
\begin{definition}
	Let $\mathcal{H}$ be a $\mathbb{Z}_{2}^{n}$-vector superspace. The complex-valued map
	$$\langle , \rangle: \mathcal{ H} \times \mathcal{ H} \rightarrow \mathbb{C}$$ is called a $\mathbb{Z}_{2}^{n}$-inner product on $\mathcal{ H}$ if the following holds:
	\begin{itemize}
		\item[i)] $\langle \lambda u+\mu v,w \rangle=\lambda \langle u,w \rangle + \mu \langle v,w \rangle$, for $  \lambda , \mu \in \mathbb{C} $ and $ u,v,w \in \mathcal{ H}_a $ where $a \in \mathbb{Z}_{2}^{n}$ (Linear in the first argument);
		\item[ii)] $\langle w,v \rangle = \mathcal{B}(a,a) \overline{\langle v,w \rangle}$ for $ v,w \in \mathcal{ H}_a $ where $ a \in \mathbb{Z}_{2}^{n}$ (Hermitian symmetric);
		\item[iii)] $ \alpha(a) \langle v,v \rangle \geq 0$, for $ v \in \mathcal{ H}_a $ where $ a \in \mathbb{Z}_{2}^{n}$ (nonnegative)
		\item[iv)] $ \langle v , w \rangle =0 $, for $ v \in \mathcal{ H}_{a}$ and $ w \in \mathcal{ H}_b$, $ a \neq b $ where $ a,b \in \mathbb{Z}_{2}^{n}. $ 
	\end{itemize}
\end{definition}
The $\mathbb{Z}_{2}^{n}$-vector superspace $\mathcal{ H}$ with a $\mathbb{Z}_{2}^{n}$-inner product $\langle , \rangle$ is called a $\mathbb{Z}_{2}^{n}$-inner product space or $\mathbb{Z}_{2}^{n}$-pre-Hilbert superspace.\\
From i) and ii) it follows that $\langle , \rangle $  is antilinear or conjugate linear in the second argument, i.e. for $  \lambda , \mu \in \mathbb{C} $ and $ u,v,w \in \mathcal{ H}_a $ where $a \in \mathbb{Z}_{2}^{n}$ we have
$$ \langle u , \lambda v + \mu w \rangle = \overline{\lambda} \langle u , v \rangle + \overline {\mu} \langle u , w \rangle.$$
Let $\langle , \rangle $ be a $\mathbb{Z}_{2}^{n}$-inner product on $\mathcal{ H}$ and let  $v,w \in \mathcal{H}_a$ where $a \in \mathbb{Z}_{2}^{n}$. An inner product (in
the ordinary sense) on $\mathcal{ H}$ is defined  as fallows
$$(v,w)=\alpha(a) \langle v,w \rangle$$
The vector space $\mathcal{ H}$ with inner product $(,)$ is a pre-Hilbert space in the usual sense. We now define a norm on $\mathcal{ H}$ as fallows
\begin{align}\label{f6}
||v||=\sqrt{(v,v)}
\end{align}
The vector space $\mathcal{H}$ with this norm is a normed linear space. 
\begin{definition}
	A $\mathbb{Z}_{2}^{n}$-Hilbert superspace is a complete $\mathbb{Z}_{2}^{n}$-pre-Hilbert superspace with norm defined in \eqref{f6}.
\end{definition}
\begin{definition}
	Let $\mathcal{ H}$ be a $\mathbb{Z}_{2}^{n}$-inner product space with product $\langle , \rangle$ and let $T: \mathcal{H} \rightarrow \mathcal{ H}$ is a $\mathbb{C}$-linear map of degree $a$ where $a \in \mathbb{Z}_{2}^{n} $. The  adjoint $T^{\dag}$ from $T$ for $v \in \mathcal{ H}_{b}$ is defined as fallows
	$$\langle v,Tw \rangle =\mathcal{B}(a,b) \langle T^{\dag}v,w \rangle $$
\end{definition}
The map $T \rightarrow T^{\dag}$ is a conjugate-linear map on $\mathrm{End}_{\mathbb{C}}(\mathcal{ H})$, i.e. for every $c \in \mathbb{C}$
$$\langle (cT)^{\dag}v,w \rangle=\overline{c} \langle T^{\dag}v,w\rangle.$$
We now define adjoint $T^*$ with respect to ordinary inner product $(,)$  on $ \mathbb{Z}_{2}^{n}$-inner product space $\mathcal{ H}$ as fallows
$$(Tv,w)=(v,T^*w) \qquad for\enspace every\enspace v,w \in \mathcal{ H} $$ 
Clearly we can show that for $\mathbb{C}$-linear map $T: \mathcal{ H} \rightarrow \mathcal{ H}$ of degree $a \in \mathbb{Z}_{2}^{n}$, we have
\begin{align}\label{e5}
&&T^*=\overline{\alpha(a)}T^{\dag}.&&
\end{align}
It fallows that for $\mathbb{C}$-linear maps $T,S: \mathcal{ H} \rightarrow \mathcal{ H}$ of degrees $a,b \in \mathbb{Z}_{2}^{n}$ respectively, we have
$$T^{\dag \dag}=T, \qquad (ST)^{\dag}=\mathcal{B}(b,a) T^{\dag} S^{\dag}.$$
\subsection{smooth unitary representations}
\begin{definition}\cite{Nee11,NKHSH}\label{e7}
	Let $(\pi, \mathcal{H})$ be a unitary representation of a Lie group $G$ on Hilbert space $\mathcal{ H}$ with a smooth exponential function $exp: Lie(G) \rightarrow G$.
	\begin{itemize}
		\item[i)]
		Let vector space $\mathcal{H}^{\infty}$ denote the space of smooth vectors as follows
		$$\mathcal{ H}^{\infty}=\{v \in \mathcal{ H},\enspace \pi_{v}: G \rightarrow \mathcal{ H}\quad is \enspace smooth,\enspace \pi_{v}(g)=\pi(g)v\}$$
		map $\pi$ is smooth if space $\mathcal{ H}^{\infty}$ is dense in $\mathcal{ H}$.
		\item[ii)]
		For every $x \in Lie(G)$ and $v \in \mathcal{ H}$, Set
		\begin{align*}
		&&\overline{\mathrm{d} \pi} (x)(v):=\lim\limits_{t \rightarrow 0} \frac{1}{t}(\pi (\mathrm{exp}(tx))v-v)&&
		\end{align*}	
		and set $\mathcal{D}_{x}$  consists of all vectors $v \in \mathcal{ H}$ such that $\frac{\mathrm{d}}{\mathrm{d}t}|_{t=0} \pi(exp_{G}(tx)v)$ exists. For the one-parameter group $ \pi (\mathrm{exp}(tx)) $, $\mathcal{D}_{x}$ is the domain of the infinitesimal generator of unbounded operator $ \overline{\mathrm{d}\pi}(x)$ and $ -i \overline{\mathrm{d}\pi}(x) $ is the self-adjoint generator. The domain of the unbounded operator $ \overline{\mathrm{d}\pi}(x)$ is defined by $\mathcal{D}(\overline{\mathrm{d}\pi}(x))$.
		\item[iii)]
		Let $D^1:=\bigcap\limits_{x \in \mathfrak{g}}D_x$ and $ D^n:=\{v \in D^1: (\forall x \in \mathfrak{g}) \overline{\mathrm{d}\pi}(x)v \in D^{n-1}\} $ for $ n > 1 $, then we put  $D^{\infty}:=\bigcap\limits_{n \in \mathbb{N}} D^n$. 
	\end{itemize}
\end{definition}
The Trotter property of G and the smoothness
of the representation imply that $\mathcal{D}^{\infty} =\mathcal{H}^{\infty}.$
\begin{definition}\cite{MMSH}\label{f1}
	Let $(G_0,\mathfrak{g}_{\mathbb{C}})$ be a $\mathbb{Z}_{2}^{n}$-Lie supergroup. A 4-tuple $(\pi,\mathcal{H},\mathscr{B},\rho^{\mathscr{B}})$ is a pre-representation of $(G_0,\mathfrak{g}_{\mathbb{C}})$ if the following holds
	\begin{itemize}
		\item[i)]
		$(\pi,\mathcal{H})$ is a smooth unitary representation of the Lie group $G_0$ on the $\mathbb{Z}_{2}^{n}$-graded
		Hilbert space $\mathcal{H}$ such that $\pi(g)$ preserves the $\mathbb{Z}_{2}^{n}$-grading for every $g \in G_0$.
		\item[ii)]
		$\mathscr{B}$ is a dense, $G_0$-invariant, $\mathbb{Z}_{2}^{n}$-graded subspace of $\mathcal{H}$ such that
		$$\mathscr{B} \subseteq \bigcap\limits_{x \in \mathfrak{g}_0} \mathcal{D}(\overline{\mathrm{d}\pi}(x)).$$
		\item[iii)]
		$\rho^{\mathscr{B}}:\mathfrak{g}_{\mathbb{C}}\rightarrow\mathrm{E}nd_{\mathbb{C}}(\mathscr{B})$ is a representation of the $\mathbb{Z}_{2}^{n}$-Lie superalgebra $ \mathfrak{g}_{\mathbb{C}} $, i.e., for $ x $ and $ y $ of degree $a,b$ respectively, we have
		$$\rho^{\mathscr{B}}([x,y])=\rho^{\mathscr{B}}(x)\rho^{\mathscr{B}}(y)-\mathcal{B}(a,b)\rho^{\mathscr{B}}(y)\rho^{\mathscr{B}}(x).$$
		\item[iv)]
		For every $x \in \mathfrak{g}_0$, $\rho^{\mathscr{B}}(x)=\overline{\mathrm{d}\pi}(x)|_{\mathscr{B}}$
		and $\rho^{\mathscr{B}}(x)$ is essentially skew adjoint, where  
		$\overline{\mathrm{d}\pi}(x)$ denotes the infinitesimal generator of the one-parameter
		group $t \mapsto \pi(e^{tx})$.
		\item[v)]
		For every $x \in \mathfrak{g}_{\mathbb{C}}$, $\rho^{\mathscr{B}}(x)^{\dag}=-\rho^{\mathscr{B}}(x)$, i.e., for x of degree a, $\overline{\alpha(a)}^{\frac{1}{2}}\rho^{\mathscr{B}}(x)$ is a symmetric operator, i.e., $\overline{\alpha(a)} \rho^{\mathscr{B}}(x) \subseteq \rho^{\mathscr{B}}(x)^*$.
		\item[vi)]
		$\rho^{\mathscr{B}}$ is a homomorphism of\enspace  $G_0$-modulus, i.e., for every $x \in \mathfrak{g}_{\mathbb{C}}$ and $g \in G_0$, $\pi(g)\rho^{\mathscr{B}}(x)\pi(g)^{-1}=\rho^{\mathscr{B}}(\mathrm{Ad}(g)x)$.
	\end{itemize}
\end{definition}
\begin{definition}\cite{MMSH}\label{f5}
	Let $(G_0,\mathfrak{g}_{\mathbb{C}})$ be a $\mathbb{Z}_{2}^{n}$-Lie supergroup. A smooth unitary representation of $(G_0,\mathfrak{g}_{\mathbb{C}})$ is a triple $(\pi,\rho^{\pi},\mathcal{H})$ with the following properties.
	\begin{itemize}
		\item[i)]
		$(\pi,\mathcal{H})$ is a smooth unitary representation of the Lie group $G_0$ on the $\mathbb{Z}_{2}^{n}$-graded
		Hilbert space $\mathcal{H}$ such that $\pi(g)$ preserves the $\mathbb{Z}_{2}^{n}$-grading for every $g \in G_0$.
		\item[ii)]
		For the space  of smooth vectors; $\mathcal{H}^{\infty}$,  $\rho^{\pi}:\mathfrak{g}_{\mathbb{C}}\rightarrow\mathrm{E}nd_{\mathbb{C}}(\mathcal{H}^{\infty})$ is a representation of the $\mathbb{Z}_{2}^{n}$-Lie superalgebra $ \mathfrak{g}_{\mathbb{C}} $, i.e., for $ x $ and $ y $ of degrees $a,b$ respectively, one has
		$$\rho^{\pi}([x,y])=\rho^{\pi}(x)\rho^{\pi}(y)-\mathcal{B}(a,b)\rho^{\pi}(y)\rho^{\pi}(x).$$
		\item[iii)]
		For every $x \in \mathfrak{g}_0$, $\rho^{\pi}(x)=\overline{\mathrm{d}\pi}(x)|_{\mathcal{H}^{\infty}}$.
		\item[\text{iv})]
		$\rho^{\pi}(x)^{\dag}=-\rho^{\pi}(x)$ For every $x \in \mathfrak{g}$, i.e., $\overline{\alpha(a)}^{\frac{1}{2}}\rho^{\pi}(x)$ is a symmetric operator, i.e.,
		$\overline{\alpha(a)} \rho^{\pi}(x) \subseteq \rho^{\pi}(x)^*$  for $x$ of degree $a$.
		\item[v)]
		$\rho^{\pi}$ is a homomorphism of\enspace  $G_0$-modulus, i.e., $\pi(g)\rho^{\pi}(x)\pi(g)^{-1}=\rho^{\pi}(\mathrm{Ad}(g)x)$ for every $x \in \mathfrak{g}_{\mathbb{C}}$ and $g \in G_0$.
	\end{itemize}
\end{definition}
\section{ Stability Theorem}\label{e12}
In this section, the stability theorem has been proved without extra condition different from what is said in \cite{MMSH}.\\
Let $a=(a_1,a_2,...,a_n)$, $b=(b_1,b_2,..,b_n)$ $\in \mathbb{Z}_{2}^{n}$, we recall the following concepts:
\begin{itemize}
	\item [i)]$\mathcal{B}(a,b)=(-1)^{\langle a,b\rangle}$
	where $ \langle a,b \rangle=a_1b_1+a_2b_2+...+a_nb_n $.
	\item[ii)]
	$\alpha(a)=e^{i \frac{\pi}{2} \textbf{u}(a)}$, where $ \textbf{u}(a)=|\{1\leq j\leq n\enspace : \enspace a_j=\bar{1}\}| $
\end{itemize}
The following auxiliary lemma is needed. See \cite[Lemma 2.5]{MSNKHSH} for a proof.
\begin{lemma}\label{r1}
	Let $ P_1 $ and $ P_2 $ be two symmetric operators on a complex Hilbert space $ \mathcal{H} $ such
	that $ \mathcal{D}(P_1) = \mathcal{D}(P_2) $. Let $  \mathscr{L} \subseteq \mathcal{D}(P_1) $ be a dense linear subspace of $\mathcal{ H} $ such that $  P_1|_{\mathscr{L}}=P_2|_{\mathscr{L}} $.
	Assume that $ P_1|_{\mathscr{L}}$ and $P_2|_{\mathscr{L}} $  are essentially self-adjoint. Then $  P_1 = P_2 $.
\end{lemma}
Matching a unique unitary representation to a pre-representation is known as a stability
theorem. A precise statement is given in the next theorem.
\begin{theorem}[\textbf{Stability  Theorem}]\label{f2}
	Let $(\pi, \mathcal{H}, \mathscr{B}, \rho^{\mathscr{B}})$ be a pre-representation of a
	$\mathbb{Z}_{2}^{n}$-Lie supergroup $(G_0, \mathfrak{g}_{\mathbb{C}})$. Then $(\pi,\rho^{\pi}, \mathcal{H})$ is a smooth unitary representation of $(G_0,\mathfrak{g}_{\mathbb{C}})$ such that $\rho^{\pi}:\mathfrak{g}_{\mathbb{C}}\rightarrow\mathrm{E}nd_{\mathbb{C}}(\mathcal{H}^{\infty})$ is a unique linear map and $\rho^{\pi}(x)|_{\mathscr{B}}=\rho^{\mathscr{B}}(x)$.
	\begin{remark}
		At first in \cite{MMSH} this theorem is proved for $\mathbb{Z}_{2}^{n}$-Lie supergroup $(G_0, \mathfrak{g}_{\mathbb{C}})$ with the following property:\\
		for every even element $a \in \mathbb{Z}_{2}^{n}$ except for $0=(0,...,0)$
		\begin{align*}
		&&\mathfrak{g}_a=\sum\limits_{b+c=a}[\mathfrak{g}_b,\mathfrak{g}_c]&&
		\end{align*}
		where $ b,c $ are odd elements in $ \mathbb{Z}_{2}^{n} $.
	\end{remark}
	\proof
		To prove the existence of $\rho^{\pi}$, we set ${\rho}^{\pi}(x)=\overline{{\rho}^{\mathscr{B}}(x)}$ for every $x \in \mathfrak{g}_{\mathbb{C}}$, where $\overline{\rho^{\mathscr{B}}(x)}$ is a smallest closed extension on $\mathcal{H}$ with 
		$D(\overline{\rho^{\mathscr{B}}(x)})=D^{\infty}$ which is stated in the Definition \ref{e7}.\\
		\textbf{Step 1.} First we show that if $v \in \mathcal{D}^{\infty}$ then $\overline{\rho^{\mathscr{B}}(x)} v$ belongs to $\mathcal{D}^{\infty}$. It is sufficient to prove that $\overline{{\rho}^{\mathscr{B}}(x)} \in \mathcal{D}^n $ for every $x \in \mathfrak{g}_{\mathbb{C}}$ of degree a, $v \in \mathcal{D}^{n+1}$ and $n \in \mathbb{N}$.
		
		We do this by induction on n. For $y \in \mathfrak{g}_0$, $w \in \mathscr{B}$ and $v \in \mathcal{D}^{2}$, we have
		\begin{align*}
		\begin{split}
		\langle \overline{\rho^\mathscr{B}(x)} v, \overline{\mathrm{d}\pi}(y)w \rangle&=\langle \overline{\rho^\mathscr{B}(x)} v, \rho^{\mathscr{B}}(y)w \rangle=\langle v, \overline{\rho^\mathscr{B}(x)}^* \rho^\mathscr{B}(y) w\rangle\\&= \langle v, \overline{\alpha(a)}\rho^\mathscr{B}(x) \rho^\mathscr{B}(y) w\rangle
		 = \alpha(a)\langle v,\rho^{\mathscr{B}}(y) \rho^{\mathscr{B}}(x)w+\rho^{\mathscr{B}}([x,y])w\rangle\\
		&	=\alpha(a) \langle v,\rho^{\mathscr{B}}(y) \rho^{\mathscr{B}}(x)w \rangle+ \alpha(a) \langle v, \rho^{\mathscr{B}}([x,y])w\rangle\\
		&=\langle \overline{\rho^\mathscr{B}(x)} \overline{\mathrm{d}\pi}(y) v,w\rangle+ \langle \overline{\rho^{\mathscr{B}}([x,y])}v,w \rangle.
		\end{split}
		\end{align*}
		It follows that $\overline{\rho^\mathscr{B}(x)} v \in \mathcal{D}((\overline{\mathrm{d}\pi}(y)|_\mathscr{B})^*)$. Since 
		$(\overline{\mathrm{d}\pi}(y)|_\mathscr{B})^*=(\rho^\mathscr{B}(y))^*=-\overline{\rho^\mathscr{B}(y)}=-\overline{d\pi}(y)$ then $\overline{\rho^\mathscr{B}(x)} v \in \mathcal{D}(\overline{\mathrm{d}\pi}(y))$.
		
		By induction hypothesis, for
		$x_1,...,x_n \in \mathfrak{g}_{0}$ and, 
		$v \in D^{n+1}$ we obtain
		\begin{flalign*}
		\begin{split}
		\langle \overline{\mathrm{d}\pi}&(x_{n-1})...\overline{\mathrm{d}\pi}(x_1)\overline{\rho^{\mathscr{B}}(x)}v, \overline{\mathrm{d}\pi}(x_n)w\rangle=\langle v, (\overline{\mathrm{d}\pi}(x_{n-1})...\overline{\mathrm{d}\pi}(x_1)\overline{\rho^{\mathscr{B}}(x)})^*\overline{\mathrm{d}\pi}(x_n)w \rangle\\&=\langle v,\overline{\rho^{\mathscr{B}}(x)}^*\overline{\mathrm{d}\pi}(x_1)^* ...\overline{\mathrm{d}\pi}(x_{n-1})^*\overline{d\pi}(x_n)w \rangle\\
		&=\langle v,(-1)^{n+1}{\overline{\alpha(a)}}\enspace\overline{\rho^{\mathscr{B}}(x)}\overline{\mathrm{d}\pi}(x_1) ...\overline{\mathrm{d}\pi}(x_{n-1})\overline{\mathrm{d}\pi}(x_n)w \rangle\\
		&=(-1)^{n+1}\alpha(a) \langle v,\rho^\mathscr{B}(x) \rho^\mathscr{B}(x_1)...\rho^\mathscr{B}(x_n)w \rangle\\
		&=(-1)^{n+1}\alpha(a) \langle v, \rho^\mathscr{B}([x,x_1]) \rho^\mathscr{B}(x_2)...\rho^\mathscr{B}(x_n)w+ \rho^\mathscr{B}(x_1)\rho^\mathscr{B}(x)\rho^\mathscr{B}(x_2)...\rho^\mathscr{B}(x_n)w\rangle\\
		&=\langle \overline{\mathrm{d}\pi}(x_n)...\overline{\mathrm{d}\pi}(x_2) \overline{\rho^{\mathscr{B}}([x,x_1])}v,w\rangle+\langle \overline{\mathrm{d}\pi}(x_n)...\overline{\mathrm{d}\pi}(x_2)\overline{\rho^{\mathscr{B}}(x)}\overline{\mathrm{d}\pi}(x_1)v,w \rangle
		\end{split}
		&
		\end{flalign*}
		A similar computation shows that
		$\overline{\mathrm{d}\pi}(x_{n-1})...\overline{\mathrm{d}\pi}(x_1)\overline{\rho^{\mathscr{B}}(x)}v \in D(\overline{\mathrm{d}\pi}(x_n))$.
		Consequently,
		$\overline{\rho^{\mathscr{B}}(x)}v \in \mathcal{D}^n$.
		
		Also the Trotter property of G and the smoothness
		of the representation imply that $\mathcal{D}^{\infty} =\mathcal{H}^{\infty} $.
		
		\textbf{Step 2.}
		We show that $\rho^{\pi}:\mathfrak{g}_{\mathbb{C}}\rightarrow\mathrm{E}nd_{\mathbb{C}}(\mathcal{H}^{\infty})$ is a representation of the $\mathbb{Z}_{2}^{n}$-Lie superalgebra $ \mathfrak{g}_{\mathbb{C}} $.
		For $x \in \mathfrak{g}_0$, $v \in \mathscr{B}$ and $c \in \mathbb{R}$, we have $$\overline{\rho^{\mathscr{B}}(cx)}v=\overline{\mathrm{d}\pi}(ax)v=c \overline{\mathrm{d}\pi}(x)v=c\overline{\rho^{\mathscr{B}}(x)}v.$$
		For $x \in \mathfrak{g}_{\mathbb{C}}$ of degree a and $v \in \mathcal{D}^{\infty}$ if we set 
		$$\mathscr{L}=\mathscr{B},\qquad P_1v=\overline{\alpha(a)}^{\frac{1}{2}}\overline{\rho^{\mathscr{B}}(cx)}v\qquad \text{and}\qquad P_2=\overline{\alpha(a)}^{\frac{1}{2}}c\overline{\rho^{\mathscr{B}}(x)}v$$
		by Lemma \ref{r1}, $P_1=P_2$ for every $  v \in \mathcal{D}^{\infty}.$
		For $x,y \in \mathfrak{g}_{\mathbb{C}}$, both of degree $ a $, and $v \in \mathcal{D}^{\infty}$, by apply Lemma \ref{r1}, a similar
		reasoning show that,
		$$\overline{\rho^{\mathscr{B}}(x+y)}v=\overline{\rho^{\mathscr{B}}(x)}v+\overline{\rho^{\mathscr{B}}(y)}v.$$
		Now, in order to prove $ \rho^{\mathscr{B}} $ preserves the Lie braket, let $x, y \in \mathfrak{g}_{\mathbb{C}}$ be of degree a, b, respictively.
		we define two operators $T_1$ and $T_2$ with
		domains $ \mathcal{D}(T_1) = \mathcal{D}(T_2) = \mathcal{D}^{\infty} $ as follows. For $  v \in \mathcal{D}^{\infty}  $ we set
		\begin{align*}
		\begin{split}
		&T_1v=(\mathcal{B}(a,b)\overline{\alpha(a)}\overline{\alpha(b)})^{\frac{1}{2}} \overline{\rho^{\mathscr{B}}([x,y])}v,\\
		&T_2v=(\mathcal{B}(a,b)\overline{\alpha(a)}\overline{\alpha(b)})^{\frac{1}{2}}(\overline{\rho^{\mathscr{B}}(x)}\enspace\overline{\rho^{\mathscr{B}}(y)} v-\mathcal{B}(a,b)\overline{\rho^{\mathscr{B}}(y)}\enspace\overline{\rho^{\mathscr{B}}(x)}v)
		\end{split} 
		&
		\end{align*}
		where let $ x, y \in \mathfrak{g}_{\mathbb{C}} $ be of degree a, b, respictively.\\
		Then $ T_1 $ and $ T_2 $ are both symmetric, $T_1|_{\mathscr{B}}=T_2|_{\mathscr{B}}$
		and since
		\begingroup\makeatletter\def\f@size{9.2}\check@mathfonts
		\begin{flalign*}
		\begin{split}
		{T_1|_{\mathscr{B}}}^*&=((\mathcal{B}(a,b)\overline{\alpha(a)}\overline{\alpha(b)})^{\frac{1}{2}}\rho^{\mathscr{B}}([x,y]))^*=\overline{\mathcal{B}(a,b)^{\frac{1}{2}}}{\alpha(a)}^{\frac{1}{2}}{\alpha(b)}^{\frac{1}{2}}\rho^{\mathscr{B}}([x,y])^*\\&=\overline{\mathcal{B}(a,b)^{\frac{1}{2}}}{\alpha(a)}^{\frac{1}{2}}{\alpha(b)}^{\frac{1}{2}}\overline{\rho^{\mathscr{B}}([x,y])^*}=-\overline{\mathcal{B}(a,b)^{\frac{1}{2}}}{\alpha(a)}^{\frac{1}{2}}{\alpha(b)}^{\frac{1}{2}}\overline{\alpha(a+b)}\overline{\rho^{\mathscr{B}}([x,y])^{\dag}}\\
		& 	=\overline{\mathcal{B}(a,b)^{\frac{1}{2}}}{\alpha(a)}^{\frac{1}{2}}{\alpha(b)}^{\frac{1}{2}}\overline{\alpha(a+b)}\overline{\rho^{\mathscr{B}}([x,y])}=\overline{\mathcal{B}(a,b)^{\frac{1}{2}}}{\alpha(a)}^{\frac{1}{2}}{\alpha(b)}^{\frac{1}{2}}\overline{\alpha(a)}\overline{\alpha(b)}\mathcal{B}(a,b)\overline{\rho^{\mathscr{B}}([x,y])}\\
		&=\overline{\alpha(a)}^{\frac{1}{2}}\alpha(a)^{\frac{1}{2}}\overline{\alpha(b)}^{\frac{1}{2}}\alpha(b)^{\frac{1}{2}}\overline{\mathcal{B}(a,b)^{\frac{1}{2}}}{\alpha(a)}^{\frac{1}{2}}{\alpha(b)}^{\frac{1}{2}}\overline{\alpha(a)}\overline{\alpha(b)}\mathcal{B}(a,b)^{\frac{1}{2}}\mathcal{B}(a,b)^{\frac{1}{2}}\overline{\rho^{\mathscr{B}}([x,y])}\\
		&=\overline{\alpha(a)}^{\frac{1}{2}}\overline{\alpha(b)}^{\frac{1}{2}}\mathcal{B}(a,b)^{\frac{1}{2}}\overline{\rho^{\mathscr{B}}([x,y])}=\overline{T_1|_{\mathscr{B}}}
		\end{split}
		&
		\end{flalign*}
		\endgroup
		the operator $ T_1|_{\mathscr{B}} $
		is essentially self-adjoint. Then Lemma \ref{r1} implies that $ T_1 = T_2 $.
		
		\textbf{Step 3.} 
		To prove 	$\rho^{\pi}$ is a homomorphism of\enspace  $G_0$-modulus, we apply Lemma \ref{r1} with
		$$\mathscr{L}=\mathscr{B},\qquad P_1=\overline{\alpha(a)^{\frac{1}{2}}}\pi(g)\overline{\rho^{\mathscr{B}}(x)} \pi(g)^{-1}, \qquad P_2=\overline{\alpha(a)}^{\frac{1}{2}}\overline{\rho^{\mathscr{B}}(Ad(g)x)}.$$
		For proving that the smooth unitary representation
		$(\pi, \rho^{\pi}, \mathcal{ H})$
		satisfies $\rho^{\pi}(x)|_{\mathscr{B}}=\rho^{\mathscr{B}}(x)$, for every $ x \in \mathfrak{g}_{\mathbb{C}} $ of degree $ a $, apply Lemma \ref{r1} with
		$ \mathscr{L}=\mathscr{B},\enspace P_1=\overline{\alpha(a)}^{\frac{1}{2}}\rho^{\pi}(x)|_{\mathcal{H}^{\infty}}\enspace \text{and}\enspace  P_2=\overline{\alpha(a)}^{\frac{1}{2}}\overline{\rho^{\mathscr{B}}(x)} $. This implies that $\rho^{\pi}(x)$ is unique.
\end{theorem}
\section{positive definite $\mathbb{Z}_{2}^{n}$-superfunctions
}\label{e9}
In this section, the concepts of uniformly boundedness and locally positive definiteness of $ \mathbb{Z}_{2}^{n} $-superfunctions have been discussed.	
\subsection{
	local superfunction on $\mathbb{Z}_{2}^{n}$-Lie supergroup
}
Let $f_{V} \in \mathrm{Hom}_{\mathfrak{g}_{0}}(\mathfrak{U}(\mathfrak{g}_{\mathbb{C}}),C_{G_{0}}^{\infty}(V))$ be a $ \mathbb{Z}_{2}^{n} $-superfunction on the $ \mathbb{Z}_{2}^{n} $-Lie supergroup $ (G_{0},\mathfrak{g}_{\mathbb{C}}) $, where $V$ is an open neighborhood of 1, identity element of Lie group $  G_0 $, and $\mathfrak{U}(\mathfrak{g}_{\mathbb{C}})$ is an enveloping algebra of $\mathbb{Z}_{2}^{n}$-complex superalgebra $\mathfrak{g}_{\mathbb{C}}=\mathfrak{g} \otimes_{\mathbb{R}}\mathbb{C}$, where $\mathfrak{g}$ is a $\mathbb{Z}_{2}^{n}$-real superalgebra. From now on we assume that $V=UU^{-1}$, where $U$ is an open neighborhood of 1 in Lie group $G_{0}$ (this can be
done after shrinking V if necessary). Indeed, let  $m:G_{0} \times G_{0}\rightarrow G_{0}$ be multiplication operation of Lie group $G_{0}$. Then there exists an open neighborhood $U \subseteq G_{0}$ such that image $U \times U$ is a subset V. consider the following maps
\begin{equation*}
\begin{array}{c}
G_0 \times G_0 \xrightarrow{1 \times r}G_0 \times G_0 \xrightarrow{m} G_0\\
(g,h) \mapsto (g,h^{-1}) \mapsto gh^{-1}
\end{array}
\end{equation*}
Since $ V $ is an open subset of $ G_0 $, 
$ (m \circ (1 \times r))^{-1}(V) $ is an open subset of $G_{0} \times G_{0}$ of $(1,1)$. Therefore there exist open neighborhoods $U_{1}$, $U_{2}$ of 1 such that $U_{1} \times U_{2}$ is subset of $ (m \circ (1 \times r))^{-1}(V) $. Set $U:= U_{1} \cap U_{2}$, then $U \times U$ is mapped into the subset V, i.e., $UU^{-1} \subseteq V$.\\
\subsection{
	uniformly bounded 	$ \mathbb{Z}_{2}^{n} $-superfunction
}
\begin{definition}\label{e6}
	A $\mathbb{Z}_{2}^{n}$-superfunction $f_{V} \in \mathrm{Hom}_{\mathfrak{g}_{0}}(\mathfrak{U}(\mathfrak{g}_{\mathbb{C}}),C_{G_{0}}^{\infty}(V))$ (or in short $f_{V} \in \mathrm{Hom}$) is called uniformly bounded if the family $\{f_{V}(D)\}_{D \in \mathfrak{U}(\mathfrak{g}_{\mathbb{C}})}$ of functions is uniformly bounded, i.e. there exists a real number $ C $ such that  $$|f_{V}(D)(x)| < C$$
	for every $x \in V$ and $D \in \mathfrak{U}(\mathfrak{g}_{\mathbb{C}})$. We denote the set of these superfunctions by $C_{ub}^{\infty}(V)$.
\end{definition}
Let $V$ be an open subset of $G_{0}$. The following is an isomorphism between superalgebras.
\begin{align*}
s:\mathrm{Hom}_{\mathfrak{g_{0}}}(\mathfrak{U}(\mathfrak{g}_{\mathbb{C}}),C^{\infty}(V)))\rightarrow C^{\infty}(V)\otimes S \mathfrak{g}_{1}^{*},\enspace f \mapsto \sum_{\stackrel {1\leq i_{1}<...<i_{k}\leq n}{\enspace 0 \leq k \leq n}}  f_{i_{1}...i_{k}}e_{i_{1}^{*}...e_{i_{k}}^{*}}	
\end{align*}
where $\{e_{1},...,e_{n}\}$ and $\{e_{i_{1}}... e_{i_{k}}\}_{i_{1} <...< i_{k}}$ are the base for $\mathfrak{g}_{1}$ and $S^{k}\mathfrak{g}_{1}$ respectively and $f(e_{i_{1}}... e_{i_{k}})=f_{i_{1}...i_{k}}$.
Note that by $\mathfrak{U}(\mathfrak{g}_{\mathbb{C}}) \cong \mathfrak{U}(\mathfrak{g_{0}})\otimes S(g_1)$ where $ g_1=\bigoplus_{a\in \mathbb{Z}_2^n-{0}} g_a $ and $ S(g_1) $ is the symmetric algebra of $ g_1 $ (see \cite{Manin1988}). For a proof in case of n=1 see \cite[Lem 7.4.8]{carmeli2011mathematical} one may consider $ e_{i_{1}}...e_{i_{k}} $ as an element of $ \mathfrak{U}(\mathfrak{g}_{\mathbb{C}}) $.\\
Suppose $h e_{i_{1}}^{*}...e_{i_{k}}^{*}$ is an element of $ C^{\infty}(V)\otimes S \mathfrak{g}_{1}^{*} $. Set
$$\tilde{h}: \mathfrak{U}(\mathfrak{g}_{\mathbb{C}}) \rightarrow C^{\infty}(V)$$
with $ \tilde{h}(e_{i_{1}}...e_{i_{k}})=h $ and $ \tilde{h}(e_{j_{1}}... e_{j_l})=0$ whenever $  \{j_1, ..., j_l\} $ is not equal to $  \{i_1, ..., i_k\} $. In addition for each $D_{0} \in \mathfrak{U}(\mathfrak{g_{0}})$ set
$$\tilde{h}(D_{0}\otimes e_{i_{1}}... e_{i_{k}})=L_{D_{0}}\tilde{h}(e_{i_{1}}...e_{i_{k}}).$$
where $  L_{D_0} $ is left-invariant differential operator corresponding to $D_0$.\\
We are looking for an interpretation of uniformly bounded functions in the topology of seminorms. At first, the topology of seminorms will be introduced and then a comparison will be made between $J$-adic topology and topology of seminorms.\\

Suppose $\chi$ is a vector space and $\rho$ is a family of seminorms on $\chi$. Let $\mathfrak{J}$ be a topology on $\chi$ generated by the subbase consisting of sets as $\{x:p(x-x_{0})<\varepsilon\}$, for $p \in \rho$, $x_{0} \in \chi$, and $\varepsilon> 0$.
Therefore, a subset $W$ of $\chi$ is open if and only if for every $x_{0}$ belongs $W$, there exists $p_{1},...,p_{n} \in \rho$ and $\varepsilon_{1},...,\varepsilon_{n}>0$ such that $\bigcap\limits_{j=1}^{n}\{x \in \chi:p_{j}(x-x_{0})<\varepsilon_{j} \} \subseteq W$. The vector space $\chi$ with this topology is a topological vector space.\cite{Conway1985}
\\
Let $ V $ be an open subset of the Lie group $ G_0 $. Set $\mathrm{A(V)}:= C^{\infty}(V)\otimes S \mathfrak{g}_{1}^{*}$. For every compact subset $Z \subseteq V$ and every differential operator $\partial=\sum D_{I}\frac{\partial}{\partial e_{I}^{*}}$ in $ \mathrm{Diff}\enspace \mathrm{A(V)} $, where ${e_{1},...,e_{n}}$ is a basis for $\mathfrak{g}_{1}$, $ I=\{i_1<...<i_k\}\subset \{1, ..., n\} $ is a multi-index and $ D_I $ is a differential operator on open set $ V $, define a seminorm $|.|_{Z,\partial}$ on $ \mathrm{A(V)} $ as follows; 
\begin{align}\label{f1}
|f|_{Z,\partial}=\sup\limits_{p \in Z}|\widetilde{\partial f}(p)|
\end{align}
Where $f=\sum f_{I}e^{*I} \in \mathrm{A(V)} $ and $e^{*I}=e^{*i_{1}}... e^{*i_{k}}$. The set of all such seminorms induce a topology on $ \mathrm{A(V)} $.\\
For example, suppose $  G_0=\mathbb{R} $ and $ V $ is an open interval in $ \mathbb{R} $ and $ x_1 $ is a coordinate function on $ V $ and $Z=\{p\} \subset V$. In this case, the seminorm $  |.|_{Z, \partial} $ is as follows:
\begin{equation*}
|f|_{Z,\partial}=|\sum f_{I}e^{*I}|_{Z,\partial}= \left\{
\begin{array}{ll}
|f_{1}(p)| & \text{if }\quad\partial=\frac{\partial}{\partial e_{1}^{*}}\\
|f_{12}(p)| & \text{if }\quad\partial=\frac{\partial}{\partial e_{1}^{*}}\frac{\partial}{\partial e_{2}^{*}}\\
|\frac{\partial f_{\varnothing}}{\partial x_{1}}(p)| & \text{if }\quad\partial=\frac{\partial}{\partial x_{1}}
\end{array} \right.
\end{equation*}
\begin{definition}\cite[Def 1.2]{Conway1985}
	A topological vector space is called a locally convex space if its topology is defined by a family of seminorms 	$\rho$ such that 	$\bigcap\limits_{p \in \rho}\{x: p(x)=0\}=(0).$ 
\end{definition}
This condition guarantees that this space is Hausdorff space \cite{Conway1985}.\\
According to the definition, a function $f$ is called bounded by a positive real number $ M $ in the topology of seminorms, denoted by $|f| \prec  M$, whenever for each differential operator $\partial$ and each compact subset $Z$
$$|f|_{Z,\partial}=\sup\limits_{p \in Z}|\widetilde{\partial f}(p)|<N_{\partial,Z}M$$
where $N_{\partial,Z}$ is a real number depends on $\partial$ and $Z$.\\
Here we give an interpretation of uniformly bounded functions in term of the topology of seminorms.
\begin{proposition}
	Let $\mathbb{Z}_{2}^{n}$-superfunction $f \in \mathrm{Hom}_{\mathfrak{g}_{0}}(\mathfrak{U}(\mathfrak{g}_{\mathbb{C}}),C_{G_{0}}^{\infty}(V))$ be uniformly bounded. Then $s(f)$ is bounded  in term of the topology of seminorms. 
\end{proposition}
\proof
	Suppose $\mathbb{Z}_{2}^{n}$-superfunction $f \in \mathrm{Hom}$ is uniformly bounded that means the family $\{f(D)\}_{D \in \mathfrak{U}(\mathfrak{g}_{\mathbb{C}})}$ of smooth functions is uniformly bounded. That is, for every $x \in V$ and $D \in \mathfrak{U}(\mathfrak{g}_{\mathbb{C}}) $ we have $|f(D)(x)|< C$, where $C$ is a real number independent of $D$ and $x$, we must show that $s(f)$ is bounded by $C$ in the topology of seminorms \eqref{f1}.
	That is
	$$|s(f)|=|\sum\limits_{1\leq i_{1}<...<i_{k}\leq n,\enspace 0 \leq k \leq n}f_{i_{1}...i_{k}}e_{i_{1}^{*}...e_{i_{k}}^{*}}| \prec C$$
	For this we show that for each compact subset $Z$ and each differential operator $\partial=\sum a_{\mu I}\frac{\partial}{\partial x_{\mu}}\frac{\partial}{\partial e_{I}^{*}}$, $ I $ and $ \mu $ are multiple indices, there exists  a real number  $ N_{\partial, Z} $ such that $$|s(f)|_{Z,\partial}=\sup\limits_{p \in Z \subset U}|\widetilde{\partial s(f)}(p)|<N_{\partial,Z}C,$$ To this end, first consider the following inequality
	\begin{align*}
	\sup\limits_{p \in Z \subset U}|\widetilde{\partial s(f)}(p)|&=\sup\limits_{p \in Z \subset U}|\sum\limits_{\mu I}a_{\mu I}\frac{\partial f_{I}}{\partial x_{\mu}}|\\& <\max\sup\{{a}_{\mu I}\}C
	\end{align*}
	It shows that it is enough to put $N_{\partial,Z}:=\max\sup\{{a}_{\mu I}\}$.
\begin{remark}
	Here, we compare the two topologies on $  A(V)  $. In fact, we show that the $J$-adic topology on $\mathrm{A(V)}$ is finer than the seminorm topology. For this purpose, it is enough to show that for every $h \in \bigcap\limits_{j=1}^{n}\{f \in \mathrm{A(V)}:p_{j}(f)<\varepsilon_{j} \}$ where $p_{j}(f)=|f|_{Z_{j},\partial_{j}}$, there exists $t$ such that $h+J^{t}\subseteq \bigcap\limits_{j=1}^{n}\{f \in \mathrm{A(V)}:p_{j}(f)<\varepsilon_{j} \}$.\\
	To this end  it is enough to show that $|h+g|_{Z_{j},\partial_{j}}<\varepsilon_{j}$ for each $g \in J^{t}$.\\
	Since $h \in \bigcap\limits_{j=1}^{n}\{f \in \mathrm{A(V)}:p_{j}(f)<\varepsilon_{j} \}$, then $|h|_{Z_{j},\partial_{j}}<\varepsilon_{j}$.\\
	Let $ k $ be the maximum odd degree of $  \partial_j $, $ 1 \leq j \leq n $. By odd degree  of a differential operator $ \partial $, we mean the highest degree of the monomials including odd derivations. Then for each $ t>k $, one has $  |g|_{Z_j, \partial_j}=0 $ if $ g\in J^t $. Thus $ |h+ g|_{Z_{j},\partial_{j}} \leq |h|_{Z_{j},\partial_{j}}+ |g|_{Z_{j},\partial_{j}}<\varepsilon_{j} $.
\end{remark}
\subsection{local positive definite 	$\mathbb{Z}_{2}^{n}$-superfunction
}
In \cite[Sec 5]{MMSH} a binary operation is defined on $  S:= G_{0} \times \mathfrak{U}(\mathfrak{g}_{\mathbb{C}})$ as follows:
$$(g_1,D_1)(g_2,D_2)=(g_1g_2,(\mathrm{Ad}(g_{2}^{-1})(D_1))(D_2))$$
where $Ad(g)(D)$ is called the adjoint representation of $g \in G_{0}$ on $D \in \mathfrak{U}(\mathfrak{g}_{\mathbb{C}})$.It can be shown that $ S $ equipped with this operation is a monoid.\\
The neutral element of $S$ is $1_{S}:=(1_{G_{0}},1_{\mathfrak{U}(\mathfrak{g}_{\mathbb{C}})})$. Conjugation of $(g,D)$ is denoted by $ (g, D)^* $ and defined by $$(g,D)^*:=(g^{-1},\mathrm{Ad}(g)(D^*)).$$
Obviously, the morphism $ (g, D)\mapsto (g, D)^* $ is an involution on $ S $.
\\
We recall that $ \mathfrak{U}(\mathfrak{g}_{\mathbb{C}}) $ is a associative $\mathbb{Z}_{2}^{n}$-superalgebra. An element $(g,D) \in S$ is called odd (respectively even) when $D$ is an odd (even) element of $ \mathfrak{U}(\mathfrak{g}_{\mathbb{C}}) $.\\
Let $ U $ be an open subset of $ G_0 $. Henceforth, by $ S_U $ we mean $ U\times \mathfrak{U}(\mathfrak{g}_{\mathbb{C}}) $.
Note that $ S_U $ is not closed under multiplication.
The following definition is based on the Definition 5.1 in \cite{MMSH} with some variations.\\
For each $\mathbb{Z}_{2}^{n}$-superfunction $f_{UU^{-1}}\in \mathrm{Hom}_{\mathfrak{g}_{0}}(\mathfrak{U}(\mathfrak{g}_{\mathbb{C}}),C_{G_{0}}^{\infty}(UU^{-1}))$, we define complex-value map $\check{f}_{UU^{-1}}:S_{UU^{-1}} \rightarrow \mathbb{C}$ as follows
$$\check{f}_{UU^{-1}}(g,D):=f_{UU^{-1}}(D)(g).$$
\begin{definition}
	A $\mathbb{Z}_{2}^{n}$-superfunction $f_{UU^{-1}} \in \mathrm{Hom}_{\mathfrak{g}_{0}}(\mathfrak{U}(\mathfrak{g}_{\mathbb{C}}),C_{G_{0}}^{\infty}(UU^{-1}))$ is positive definite on $U$ if the following holds:
	\begin{itemize}
		\item [i)]
		The map $\check{f}_{UU^{-1}}$ is zero everywhere except $S_{UU^{-1}}^{0}=\{(g,D)\in S_{UU^{-1}},\tilde{D}=0\}$, where $\tilde{D}$ is degree $D$.
		Note that $  S_U $ and $  S_U^{-1} $ are subsets of $ S_{UU^{-1}} $. Thus for each $ s_1 $ and $ s_2 $ in $ S_U $, one has $ s_1s_2^*\in  S_{UU^{-1}} $.
		\item [ii)] For every $n\geq 1$, $c_1,...,c_n \in \mathbb{C}$ and $s_1,...,s_n \in S_{U}$, we have $$\sum\limits_{1\leq i,j \leq n} \overline{c_{i}}c_{j} \check{f}_{UU^{-1}}(s_{i}^*s_j)\geq 0.$$
	\end{itemize}
\end{definition}
\section{Constructing a pre-representation of 	$\mathbb{Z}_{2}^{n}$-Harish-Chandra pairs $ (G_0,\mathfrak{g}_{\mathbb{C}}) $.}\label{e10}
Let $f_{UU^{-1}} \in \mathrm{Hom}_{\mathfrak{g}_{0}}(\mathfrak{U}(\mathfrak{g}_{\mathbb{C}}),C_{G_{0}}^{\infty}(UU^{-1}))$ be a local positive definite $\mathbb{Z}_{2}^{n}$-superfunction which belongs to category $C_{ub}^{\infty}(UU^{-1})$.
Here a pre-representation, associated to $ f_{UU^{-1}} $,  is constructed by four steps I to IV as follows;\\
\textbf{I. Hilbert space \boldmath{$\mathcal{H}_{\check{f},U}$}}.\label{e3}
Following \cite[Page 11]{MMSH} by $\mathfrak{D}_{\check{f},U}$ we mean the vector space spanned by vectors $  K_{s, U} $, $  s\in S_U $, where
$$K_{s,U}:S_{U}\rightarrow \mathbb{C}, \qquad K_{s,U}(t):=K_{U}(t,s)=\check{f}_{UU^{-1}}(ts^*)$$
where $K_{U}:S_{U} \times S_{U} \rightarrow \mathbb{C}$. Notice that $\mathfrak{D}_{\check{f},U}$ is a $\mathbb{Z}_{2}^{n}$-graded vector space of complex-valued functions on $S_{U}$, where homogeneous parts of $\mathbb{Z}_{2}^{n}$-grading is defined as follows
$$\mathfrak{D}_{\check{f},U,a}:=\{h_{U} \in \mathfrak{D}_{\check{f},U}: h_{U}(s)=0 \qquad unless \qquad s \in S_{U,a}  \}$$
The space $\mathcal{D}_{\check{f},U}$ can be equipped with a semi-linear form as follows
\begin{align*}
&&(K_{s,U},K_{t,U}):=K_U(t,s)=K_{s,U}(t)= \check{f}_{UU^{-1}}(ts^*)&&
\end{align*}
For
$K_{s,U},K_{t,U} \in \mathfrak{D}_{\check{f},U}$, we show that
$(K_{t,U},K_{s,U})=\overline{(K_{s,U},K_{t,U})}$.
At first we show that
$K_{U}(t,s)=\overline{K_{U}(s,t)}$.\\
As $\check{f}_{UU^{-1}}$ is positive definite, then $K_{U}$ is a positive definite kernel, thus  for every $m \in \mathbb{N}$ and $s_1,...,s_m \in S_{U}$, the matrix $[K_{U}(s_i,s_j)]_{i,j=1}^{m}$ is positive definite, i.e. for every complex numbers $c_1,...,c_m \in \mathbb{C}$, one has $\sum\limits_{i,j=1}^{m}K_{U}(s_i,s_j)c_i\bar{c_j}\geq 0$.\\
we recall that if $A=[a_{ij}]$ is a positive definite matrix, then  $a_{ij}=\overline{a_{ji}}$ and eigenvalues of $ A $ are nonnegative real numbers.
\\
Let $A$ be a positive definite matrix, then for $X \neq 0$, $X^*AX>0$
where $X^*$ is conjugate transpose of $X$. If we take conjugate transpose of $X^*AX$, then $X^*A^*X>0$. Therefore $X^*AX=X^*A^*X$. In this case, we can conclude for every $X$, we have $X^*(A-A^*)X=0$.\\
If we have this relation for every $X$, then we have it for $X+Y$.
\begin{align*}
\begin{split}
(X+Y)^*(A-A^*)(X+Y)=0 &\Rightarrow X^*(A-A^*)Y+Y^*(A-A^*)X=0\\
& \Rightarrow X^*(A-A^*)Y-(X^*(A-A^*)Y)^*=0\\
&\Rightarrow X^*(A-A^*)Y-\overline{(X^*(A-A^*)Y)}=0
\end{split}
&
\end{align*}
Then $X^*(A-A^*)Y \in \mathbb{R}$. we shall show that for every $X$ and $Y$, $X^*(A-A^*)Y=0$.\\
Let $ \{e_i\} $ be the standard basis of $  \mathbb{C}^n $ and let $X=[e_i]$ and $Y=[e_j]$, then $ X^*(A-A^*)Y=a_{ij} $. Thus each entry of matrix $A-A^*$ is real. Also if we set $X=[e_i]$ and $Y=[\sqrt{-1}e_j]$, then one may conclude that entries $\sqrt{-1}a_{ij}$ is real. Thus $a_{ij}$ can only be zero. Therefore $A-A^*=0$. consequently $A$ is a Hermitian matrix.\\
Since $K_{U}:S_{U} \times S_{U} \rightarrow \mathbb{C}$ is a positive definite kernel, then $K_{U}(s,s)\geq 0$ for every $s \in S_{U}$ and also for every $s,t \in S_{U}$, the following matrix is Hermitian with positive determinant.
\begin{equation*}
\begin{array}{cc}
\left(
\begin{array}{cc}
K_{U}(s,s) & K_{U}(s,t)\\
K_{U}(t,s) & K_{U}(t,t)\\
\end{array} \right)
\end{array}
\end{equation*}
Then $K_{U}(t,s)=\overline{K_{U}(s,t)}$. Therefore
\begin{align*}
&&(K_{s,U},K_{t,U}):=K_U(t,s)=\overline{K_{U}(s,t)}=\overline{(K_{t,U},K_{s,U})}.&&
\end{align*}
\begin{proposition}
	For $c \in \mathbb{C}$ and $K_{s,U},K_{t,U} \in \mathfrak{D}_{\check{f},U}$, we show that
	\begin{itemize}
		\item[i)] $(K_{s,U},cK_{t,U})=\bar{c}(K_{s,U},K_{t,U})$
		\item[ii)]
		$(cK_{s,U},K_{t,U})={c}(K_{s,U},K_{t,U})$
		\item[iii)]
		($ K_{s,U}+K_{t,U},K_{m,U})=(K_{s,U},K_{m,U})+(K_{t,U},K_{m,U}) $	for $ K_{m,U} \in \mathfrak{D}_{\check{f},U}$
	\end{itemize}
\end{proposition} 
\proof
		i) Let $l \in S_{U}$, then one has
		$$ cK_{t,U}(l)=c\check{f}_{UU^{-1}}(lt^*)=\check{f}_{UU^{-1}}(clt^*)=K_{\bar{c}t,U}(l) $$
		Thus
		$$ (K_{s,U},cK_{t,U})=(K_{s,U},K_{\bar{c}t,U})=\check{f}_{UU^{-1}}(\bar{c}ts^*)=\bar{c}\check{f}_{UU^{-1}}(ts^*)=\bar{c}(K_{s,U},K_{t,U}) $$
		ii) 
		\begin{align*}
		(cK_{s,U},K_{t,U})&=(K_{\bar{c}s,U},K_{t,U})=\check{f}_{UU^{-1}}(t(\bar{c}s)^*)\\&=\check{f}_{UU^{-1}}(tcs^*)=c\check{f}_{UU^{-1}}(ts^*)=c(K_{s,U},K_{t,U})
		\end{align*} 
		iii)
		\begin{align*}
		(K_{s,U}+K_{t,U},K_{m,U})=(K_{s,U}+K_{t,U})(m)&=K_{s,U}(m)+K_{t,U}(m)\\&=(K_{s,U},K_{m,U})+(K_{t,U},K_{m,U}).
		\end{align*}	

In the last proposition it was seen that for $ h_U \in \mathfrak{D}_{\check{f},U}$ and $s \in S_{U}$, one has $h_U(s)=(h_U,K_{s,U})$. \\
We denote the completion of the pre-Hilbert space $\mathfrak{D}_{\check{f},U}$ with respect to the metric arising from Hermition inner product $ ( , ) $ by $\mathcal{ H}_{\check{f},U}$.\\\\
\textbf{II. *-representation of \boldmath{$\mathfrak{g}_{\mathbb{C}}$}}\\\\
\textbf{
	Derivative of
	\boldmath{$K_{s,U}$}
	.}
For
$ x \in \mathfrak{g}_0$,
set
\begin{align}\label{n11}
&&L_xK_{s,U}(u):= \frac{d}{dt}K_{s,U}(u(e^{tx},1))|_{t=0}&&
\end{align}
Then by definition of 
$K_{s,U}$
we have
\begin{align*}
&&L_xK_s(u):= \frac{d}{dt}\check{f}_{UU^{-1}} (u(e^{tx},1)s^*)|_{t=0}.&&
\end{align*}
Let 
$u=(p_1,D_1)$
and 
$s^*=(p_2,D_2)$
then
\begin{align*}
&&
u(e^{tx},1)s^*=(p_1,D_1)(e^{tx},1)(p_2,D_2)=(p_1e^{tx}p_2,(p_2^{-1}e^{-tx}.D_1)D_2)&
&
\end{align*}
Now we calculate the right-hand side of the second equality. For
$Y \in \mathfrak{g}$, we have
$exp \circ Ad_{P_{2}^{-1}}(Y)=c_{p_{2}^{-1}}\circ exp(Y)$
where 
$c_{p_{2}}(g)=p_2gp_{2}^{-1}$
for every 
$g \in G$. Therefore
\begin{align}\label{p12}
&&exp(Ad(p_{2}^{-1})(Y))=p_{2}^{-1}exp(Y)p_{2}&&
\end{align}
For
$Y=tx$,
by affecting $p_2$ from the left on two side of the equality
\eqref{p12}
we have
$$p_2 exp(t Ad(p_2^{-1})(x))=p_2 exp(Ad(p_{2}^{-1})(tx))=exp(tx)p_2$$
For
$Y=-tx$,
by affecting $p_{2}^{-1}$ from the right on two side of the equality
\eqref{p12}
we have
$$exp(-t Ad(p_2^{-1})(x))p_{2}^{-1}= exp(Ad(p_{2}^{-1})(-tx))p_{2}^{-1}=p_{2}^{-1}exp(-tx)$$
Therefore
\begin{align*}
u(e^{tx},1)s^*=(p_1,D_1)(e^{tx},1)(p_2,D_2)&=(p_1e^{tx}p_2,(p_2^{-1}e^{-tx}.D_1)D_2)\\
&=(p_1p_2e^{t(p_2^{-1}.x)},(e^{-t(p_2^{-1}.x)}(p_2^{-1}.D_1))D_2).
\end{align*}
Now by chain rule we obtain
\begin{flalign*}
\begin{split}
\dfrac{\mathrm{d}}{\mathrm{dt}} \check{f}_{UU^{-1}}(pe^{ty},(e^{-ty}.D_1)D_2)|_{t=0}&=\dfrac{\mathrm{d}}{\mathrm{dt}} f_{UU^{-1}}((e^{-ty}.D_1)D_2)(pe^{ty})|_{t=0}\\
&=f_{UU^{-1}}([-y,D_1]D_2)(p)+\mathrm{L}_y f_{UU^{-1}}(D_1D_2)(p)\\
&=f_{UU^{-1}}(-yD_1D_2+D_1yD_2+yD_1D_2)(p)\\
&=f_{UU^{-1}}(D_1yD_2)(p)=\check{f}(p,D_1yD_2)
\end{split}
&
\end{flalign*}
These calculations show that
\begin{flalign*}
\begin{split}
\dfrac{\mathrm{d}}{\mathrm{dt}} \check{f}_{UU^{-1}}(u(e^{tx},1)s^*)|_{t=0}&=\check{f}_{UU^{-1}}(p_1p_2,(p_2^{-1}.(D_1x))D_2)\\
&=\check{f}_{UU^{-1}}((p_1,D_1)(1,x)(p_2,D_2))=\check{f}_{UU^{-1}}(u(1,x)s^*)
\end{split}
&
\end{flalign*}
Thus
$\mathrm{L}_xK_{s,U}(u)$ exists
and indeed
$\mathrm{L}_x K_{s,U}(u)=K_{s(1,x^*),U}(u)$.
\\
Now we define a representation of $\mathbb{Z}_{2}^{n}$-Lie superalgebra $\mathfrak{g}_{\mathbb{C}}$ as follows:
\begin{align}\label{f3}
\rho:\mathfrak{g}_{\mathbb{C}}\rightarrow \mathrm{End}(\mathfrak{D}_{\check{f},U}),\qquad \rho(x)K_{s,U}:=K_{s(1,x^*),U}
\end{align}
Now we show that 
$ \rho $
is a homomorphism 
$\mathbb{Z}_{2}^{n}$-Lie superalgebra
$\mathfrak{g}_{\mathbb{C}}$, i.e.
\begin{align*}
&&\rho([x,y])=[\rho(x),\rho(y)]&&
\end{align*}
Let
$ s, t \in \mathcal{S}_U$ 
and
let $x,y$
be homogeneous elements in
$\mathfrak{g}_{\mathbb{C}}$ of degrees  
$a,b$ respectively.
Then we have
\begin{align}\label{n1}
\begin{split}
([\rho(x),\rho(y)]K_{s,U},K_{t,U})&=((\rho(x) \rho(y)) K_{s,U} - \mathcal{B}(a,b)(\rho(y) \rho(x))K_{s,U},K_{t,U})\\
&=(\rho(x)K_{s(1,y^*),U}-\mathcal{B}(a,b)\rho(y) K_{s(1,x^*),U}, K_{t,U})\\
&=(K_{s(1,y^*)(1,x^*),U}-\mathcal{B}(a,b) K_{s(1,x^*)(1,y^*),U}, K_t)\\
&=\check{f}_{UU^{-1}}(t(1,xy)s^*)-\mathcal{B}(a,b) \check{f}_{UU^{-1}}(t(1,yx)s^*)
\end{split}
&
\end{align}
For calculating the right side of the above relation, there are two points to note.\\
\textbf{The first point}\label{n2}:
Suppose $ (g,D_1),(g,D_2) \in \mathcal{S}_U$ are two arbitrary elements with the same first components, Then
\begin{align*}
\begin{split}
\check{f}_{UU^{-1}}(g,D_1)-\check{f}_{UU^{-1}}(g,D_2)&=f_{UU^{-1}}(D_1)(g)-f_{UU^{-1}}(D_2)(g)\\&=(f_{UU^{-1}}(D_1)-f_{UU^{-1}}(D_2))(g)\\&=f_{UU^{-1}}(D_1-D_2)(g)=\check{f}_{UU^{-1}}(g,D_1-D_2).
\end{split}
&
\end{align*}
\textbf{The Second point}:
The first point implies that
\begin{align*}
&&t(1,xy-\mathcal{B}(a,b)yx)s^*=t(1,xy)s^*-\mathcal{B}(a,b)t(1,yx)s^*.&&
\end{align*}
Thus we have 
\begin{align}\label{n3}
\begin{split}
\check{f}_{UU^{-1}}(t(1,xy)s^*)-\mathcal{B}(a,b)\check{f}_{UU^{-1}}(t(1,yx)s^*)&=\check{f}_{UU^{-1}}(t(1,(xy-\mathcal{B}(a,b)yx))s^*)\\&= f_{UU^{-1}}(t(1,[x, y])s^*) 
\end{split}
\end{align}
For the second equality one has to note that in
$\mathbb{Z}_{2}^{n}$-Lie superalgebra 
$\mathfrak{U}(\mathfrak{g}_{\mathbb{C}})$
for every
$x , y \in \mathfrak{g}_{\mathbb{C}}$, we have
\begin{align*}
&&[x,y]=x\otimes y -\mathcal{B}(a,b)y\otimes x&&
\end{align*}
If \eqref{n1} replaced by \eqref{n3}, then we get
\begin{align*}
\begin{split}
([\rho(x),\rho(y)]K_{s,U},K_{t,U})&=\check{f}_{UU^{-1}}(t(1,[x,y])s^*)\\&=(K_{s(1,[x,y]^*),U},K_{t,U})=(\rho([x,y])K_{s,U},K_{t,U})
\end{split}
\end{align*}
$ \rho $ is a *-representation. For this one has  to show that
$ 	(\rho(x) K_{s,U},K_{t,U})=(K_{s,U},\rho (x^*)K_{t,U}) $. To prove the equality one may start from the left side as follows
\begin{flalign*}
\begin{split}
(\rho(x) 
K_{s,U},K_{t,U})&=(K_{s(1,x^*),U},K_{t,U})=\check{f}_{UU^{-1}}(s(1,x^*)t^*)=\check{f}_{UU^{-1}}(s(t(1,x))^*)\\&=(K_{s,U},K_{t(1,x),U})=(K_{s,U},K_{t(1,(x^*)^*),U})=(K_{s,U},\rho(x^*)K_{t,U}).
\end{split}
&
\end{flalign*}
\textbf{III. Integration of \boldmath{$\rho$}}.
We show $\rho$ is integrable. By \cite[Page 5]{Nee11} it is sufficient to show that
$$\sum\limits_{r=0}^{\infty}\frac{||\rho(x_i)^{r}K_{s,U}||}{r!}< \infty$$
where $ \{x_i\} $ is a basis of the vector space $ \mathfrak{g}_{\mathbb{C}} $.\\
To this end, one has to note that
\begin{align*}
\|\rho(x_{i})^rK_{s,U}\|&=\\
&=\|K_{s(1,x_{i}^*)^r,U}\|\\
&=\|K_{s(1,x_{i}^*...x_{i}^*),U}\|\\
&= (K_{s_(1,x_{i}^*...x_{i}^*),U},K_{s(1,x_{i}^*...x_{i}^*),U})^{\frac{1}{2}}\\
&= \check{f}_{UU^{-1}}\left(s(1,x_{i}^*...x_{i}^*)(1,x_{i}...x_{i})s^*\right)^{\frac{1}{2}}\\
&=\check{f}_{UU^{-1}}\left(s(1,{x_{i}^{*}}^{r}x_{i}^r)s^*\right)^{\frac{1}{2}}
\end{align*}
\begin{align*}
s(1,{{x_{i}}^{*}}^{r}x_{i}^r)s^*&=(p,D)(1,{{x_{i}}^{*}}^{r}x_{i}^r)(p,D)^*=(p,D{{x_{i}}^{*}}^{r}x_{i}^r)(p^{-1},\mathrm{Ad}(p)(D^{*}))\\
&=(pp^{-1},\mathrm{Ad}(p)(D{{x_{i}}^{*}}^{r}x_{i}^r)\mathrm{Ad}(p)(D^{*}))=(1,\mathrm{Ad}(p)(D{{x_{i}}^{*}}^{r}x_{i}^rD^{*}))
\end{align*}
$$ \|\rho(x_{i})^rK_{s,U}\|=\check{f}_{UU^{-1}}\left(1,\mathrm{Ad}(p)(D{{x_{i}}^{*}}^{r}x_{i}^rD^{*})\right)^{\frac{1}{2}}=\left(f_{UU^{-1}}\left(\mathrm{Ad}(p)(D{{x_{i}}^{*}}^{r}x_{i}^rD^{*})\right)(1) \right)^{\frac{1}{2}}$$
Since 
$f_{UU^{-1}}$
belongs to category 
$C_{ub}^{\infty}(V)$, then for a basis
$x_1,...,x_k$
of
$\mathfrak{g}_{
	\mathbb{C}}$
it is obvious that
$$\sum\limits_{r=0}^{\infty}\frac{\|\rho(x_i)^{r}K_{s,U}\|}{r!}<\infty$$
Therefore
$\rho$ is integrable.\\\\
\textbf{IV. Representation of \boldmath{$  (G_0, \mathfrak{g}_{\mathbb{C}}) $}}.
At first we define $\pi(exph)$ for $  h\in \mathfrak{g}_{0} $ in the following way.\\
According to \eqref{f3} $\rho$ is a representation of $\mathbb{Z}_{2}^{n}$-Lie superalgebra. Since $\mathfrak{g}_{\mathbb{C}}$ is a finite dimensional, $\mathfrak{D}_{\check{f},U}$ is an equianalytic subspace of $\mathcal{H}_{\check{f},U}$, and
by the last statement of Step I on Page 15, $\mathfrak{D}_{\check{f},U}$ is also a dense subspace of  $\mathcal{H}_{\check{f},U}$.
\\
Now set $\tilde{\rho}:\mathfrak{g}_{0} \rightarrow \mathrm{U}(\mathcal{H}_{\check{f},U})$ with $\tilde{\rho}(h)=e^{\overline{\rho(h)}}$  for every $h \in \mathfrak{g}_{0}$. We show that the $ \tilde{\rho}(h) $'s are unitary operators.\\
In fact
$$\rho(h)=-\rho(h)^* \Rightarrow i\rho(h)=-i\rho(h)^* \Rightarrow (i\rho(h))^*=-i\rho(h)^*=i\rho(h)$$
then $i\rho(h)$ is self-adjoint, hence $i\rho(h)$ is symmetric and by Nelson's Theorem \cite[Th 4.8]{MMSH}, $i\rho(h)$ is essentially self-adjoint. Therefore  $\rho(h)$ is essentially skew-adjoint and this shows that $ \tilde{\rho}(h) $, for each $ h\in \mathfrak{g}_0 $, is a unitary operator. We show that $\overline{\rho(h)}$ generates an unitary one-parameter group $e^{t \overline{\rho(h)}}$ for every $t \in \mathbb{R}$. \cite[ThVIII.7]{RMSB}\\
Let $A:=\overline{\rho(h)}$, we show that $e^{tA}$ is an one-parameter group, i.e. $f: \mathbb{R}\rightarrow \mathrm{U}(\mathcal{H}_{\check{f},U})$, $f(t)=e^{tA}$ is a group homomorphism.\\
In fact by Baker-Campbell-Hausdorff formula, for every $t_1,t_2 \in \mathbb{R}$, one has 
$$e^{t_{1}A}e^{t_{2}A}=e^{t_{1}A+t_{2}A+\frac{1}{2}[t_{1}A,t_{2}A]+...}$$
since $h$ is even, so $\rho(h)$ is even. Therefore
$$[t_{1}A,t_{2}A]=t_{2}At_{1}A-t_{1}At_{2}A=0.$$
Thus f is a group homomorphism.\\
Now $\pi(exph)$ can be defined as follows:
$$\pi(exph)=\tilde{\rho}(h)$$
Assume $G_0$ is connected. A map $\pi:G_0 \rightarrow \mathrm{Aut}(\mathcal{ H}_{\check{f},U})$ is defined as follows:\\
Let $p \in G_{0}$. Hence there exists a natural number $k$ and tangent vectors $h_1,...,h_k \in \mathfrak{g}_{0}$, so that $p=exph_{1}...exph_{k}$
\cite[Corollary 3.47]{HB2015}. We define
\begin{align}\label{f4}
\pi(p):=\pi(exph_{1})...\pi(exph_{k}).
\end{align}
We show that $\pi(pq)=\pi(p)\pi(q)$ for $p,q \in G_{0}$.\\
Let $p,q \in G_{0}$ hence there exist tangent vectors $h_{1},...,h_{k} \in \mathfrak{g}_{0}$ and $l_{1},...,l_{k} \in \mathfrak{g}_{0}$ such that $p=exph_{1}...exph_{k}$ and $q=expl_{1}...expl_{k}$. By \eqref{f4} we have
$$\pi(p):=\pi(exph_{1})...\pi(exph_{k}),\qquad \pi(q):=\pi(expl_{1})...\pi(expl_{k})$$
Since
$pq=exph_1...exph_kexpl_1...expl_k$
then
\begin{align*}
&&\pi(pq)=\pi(exph_1)...\pi(exph_k)\pi(expl_1)...\pi(expl_k)&&
\end{align*}
Therefore
$\pi(pq)=\pi(p)\pi(q)$.\\
Now we show that, for
$p \in G_0$,
$\pi(p)\pi(p)^*=Id=\pi(p)^*\pi(p)$.\\
For
$p \in G_0$
there exist tangent vectors 
$h_1,...,h_k \in \mathfrak{g_{{0}}}$
such that 
$p=exph_1...exph_k$.
Consequently
\begin{align*}
\pi(p)\pi(p)^*&=\pi(exph_1)...\pi(exph_k)\left(\pi(exph_1)...\pi(exph_k)\right)^*\\&
=\pi(exph_1)...\pi(exph_k)\pi(exph_k)^*...\pi(exph_1)^*\\&=\pi(exph_1)...\pi(exph_k)\pi(-exph_k)...\pi(-exph_1)=Id.
\end{align*}
Here we show that 4-tuple 
$(\pi,\mathcal{H}_{\check{f},U},\mathfrak{D}_{\check{f},U},\rho)$
is a pre-representation of a
$\mathbb{Z}_{2}^{n}$-Harish-Chandra pair
$ 	(G_0,\mathfrak{g}_{\mathbb{C}}) $. For this we show that the conditions (i) - (vi) of Definition \ref{f1} are satisfied.\\
i) We show that
$(\pi,\mathcal{H}_{\check{f},U})$
is a smooth unitary representation of Lie group 
$G_0$
on 
the $\mathbb{Z}_{2}^{n}$-graded Hilbert space
$\mathcal{H}_{\check{f},U}$
such that 
$\pi (g)$ preserves the $\mathbb{Z}_{2}^{n}$-grading for every
$g \in G_0$.
\\
by  \cite[Def 5.1]{NKHSH}, \cite[Th 7.2]{NKH} It is sufficient to prove, for every 
$s \in S_{U}$, the map
$g \mapsto (\pi(g)K_{s,U},K_{s,U})$
from
$G_{0}$
to
$\mathbb{C}$ is smooth.
\\
Set 
$s=(g_0,D_0)$,
for every
$g \in G_0$ one has
\begingroup\makeatletter\def\f@size{9.5}\check@mathfonts
\begin{flalign*}
\begin{split}
\left(\pi(g)K_{s,U},K_{s,U}\right)&=
\left(\pi(exph_1)...\pi(exph_k)K_{s,U},K_{s,U}\right)=\left(\tilde{\rho}(h_1)...\tilde{\rho}(h_k)K_{s,U},K_{s,U}\right)\\
&=\left(e^{\overline{\rho(h_1)}}...e^{\overline{\rho(h_k)}}K_{s,U},K_{s,U}\right)=\left(\sum\limits_{r_1=0}^{\infty}\frac{\rho(h_1)^{r_1}}{r_1!}...\sum\limits_{r_1=0}^{\infty}\frac{\rho(h_k)^{r_k}}{r_k!}K_{s,U},K_{s,U}\right)\\
&=\left(\sum\limits_{r_1=0}^{\infty}...\sum\limits_{r_k=0}^{\infty}\frac{1}{r_1!}...\frac{1}{r_k!}K_{s(1,h_{k}^{*})^{r_k}...(1,h_{1}^{*})^{r_1},U},K_{s,U}\right)\\
&=\sum\limits_{r_1=0}^{\infty}...\sum\limits_{r_k=0}^{\infty}\frac{1}{r_1!}...\frac{1}{r_k!}\check{f}_{UU^{-1}}(1,Ad(g_0)(D_0h_{1}^{r_1}...h_{k}^{r_k}D_{0}^{*}))\\
&=\sum\limits_{r_1=0}^{\infty}...\sum\limits_{r_k=0}^{\infty}\frac{1}{r_1!}...\frac{1}{r_k!}f_{UU^{-1}}(Ad(g_0)(D_0h_{1}^{r_1}...h_{k}^{r_k}D_{0}^{*}))(1)
\end{split}
&
\end{flalign*}
\endgroup
Since 
$f_{UU^{-1}}$
is uniformly bounded, then this series is convergent. Therefore map
$g \mapsto (\pi(g)K_{s,U},K_{s,U})$
from
$G_{0}$
to
$\mathbb{C}$ is smooth.
In order to prove forth condition of the definition of pre-representation , we have to show that
$\rho(x)=\overline{\mathrm{d}\pi}(x)\arrowvert_{\mathfrak{D}_{\check{f},U}}$
and
$\rho(x)$ essentially skew adjoint for
$x \in \mathfrak{g}_{0}$
\\
For the first part, let
$K_{s,U} \in \mathfrak{D}_{\check{f},U}$
\begin{align*}
\begin{split}
\overline{d\pi}(x)K_{s,U}&=\lim\limits_{t \rightarrow 0}\frac{1}{t}(\pi(e^{tx})K_{s,U}-K_{s,U})=\frac{d}{dt}|_{t=0}\pi(e^{tx})K_{s,U}\\&=\frac{d}{dt}|_{t=0}Exp(\rho(tx))K_{s,U}=Exp'_{0}\frac{d}{dt}\rho(tx)K_{s,U}=\rho(x)K_{s,U}
\end{split}
&
\end{align*}
In order to prove the second part, it is first necessary to prove the fifth condition of the definition of pre-representation.\\
To this end one has to show
($\rho(x)^{\dag}=-\rho(x)$). By relations 
$\rho(x)^*=\overline{\alpha(|\rho(x)|)}\rho(x)^{\dag}$, 
$ \rho(x)^*K_{s,U}=K_{s(1,x),U} $
and
$ K_{s(1,cx^*),U}=\overline{c}K_{s(1,x^*),U} $, 
for
$c \in \mathbb{C}$ one has
\begin{align*}
\rho(x)^{\dag}K_{s,U}&=\alpha(a)\rho(x)^*K_{s,U}=\alpha(a)K_{s(1,x),U}=\alpha(a)K_{s(1,-{\alpha(a)}x^*),U}\\
&=-\alpha(a)\overline{\alpha(a)}K_{s(1,x^*),U}=-K_{s(1,x^*),U}
=-\rho(x)K_{s,U}.
\end{align*}
In order to prove the second part, we use a densely defined operator 
$\rho(x)$
on a Hilbert space 
$ \mathcal{H}_{\check{f},U} $,  
the fifth condition of definition of pre-representation 
($\rho(x)^{\dag}=-\rho(x) $)
and a relation expressed in \eqref{e5}
$\rho(x)^*=\overline{\alpha(|\rho(x)|)}\rho(x)^{\dag} $
for
$x$ of degree zero.
\begin{align*}
&&\overline{\rho(x)}=\overline{-\rho(x)^{\dag}}=\overline{-\rho(x)^*}=-\overline{\rho(x)^*}=-\rho(x)^*.&&
\end{align*}
We show that
$\rho(x)^{\dag}=-\rho(x)$ for  
$x \in \mathfrak{g}$. By 
\eqref{e4}
and
\eqref{e5}
for 
$K_{s,U} \in \mathfrak{D}_{\check{f},U}$
we have
\begin{align*}
\begin{split}
\rho(x)^{\dag}K_{s,U}&=\alpha(a)\rho(x)^*K_{s,U}=\alpha(a)K_{s(1,x),U}=\alpha(a)K_{s(1,-{\alpha(a)}x^*),U}\\
&=-\alpha(a)\overline{\alpha(a)}K_{s(1,x^*),U}=-K_{s(1,x^*),U}
=-\rho(x)K_{s,U}.
\end{split}
&
\end{align*}
We show that the sixth condition of the definition of pre-representation is established 
$\left( \pi(g) \rho(x) \pi(g)^{-1}=\rho(\mathrm{Ad}(g)x)\right)$
for
$ g \in G_0$
and
$x \in \mathfrak{g}_{\mathbb{C}}$.
\\
For 
$g \in G_0$
there exist tangent vectors  
$h_1,...,h_n \in \mathfrak{g_{0}}$
such that
$g=exph_1...exph_n$.
Then for
$g^{-1} \in G_0$,
we can write 
$g^{-1}=exp(-h_n)...exp(-h_1)$. We prove this by induction on $ n $.
For
$K_{s,U} \in \mathfrak{D}_{\check{f},U}$
we show that
$ \pi(exp h)\rho(x)\\\pi(exp h)^{-1}=\rho\left(Ad(exp h)x\right) $.
Consider the commutative diagrams
\begin{equation*}
\xymatrix{
	G_0 \ar[r]^{\mathrm{Ad}} &
	GL(\mathfrak{g_{0}})  \\
	\mathfrak{g_{{0}}} \ar[u]^{exp}\ar[r]^{\mathrm{ad}} & \mathrm{End}(\mathfrak{g_{0}}) \ar[u]^{\mathrm{EXP}} }\quad\quad
\xymatrix{
	GL(\mathfrak{D}_{\check{f},U}) \ar[r]^{\mathrm{Ad}_{\mathfrak{D}}} &
	GL(\mathrm{End}(\mathfrak{D}_{\check{f},U}))  \\
	\mathrm{End}(\mathfrak{D}_{\check{f},U}) \ar[u]^{exp_{\mathfrak{D}}}\ar[r]^{\mathrm{ad}_{\mathfrak{D}}} & \mathrm{End}(\mathrm{End}(\mathfrak{D}_{\check{f},U})) \ar[u]^{\mathrm{EXP}_{\mathfrak{D}}}
}
\end{equation*}
We have
\begin{align*}
\begin{split}
\rho(\mathrm{Ad}(exph)x)&=\rho(\mathrm{Exp}(\mathrm{ad}h)x)=\rho\left(\sum\limits_{n=0}^{\infty}\frac{1}{n!}(\mathrm{ad}h)^{n}x\right)=\sum\limits_{n=0}^{\infty}\frac{1}{n!}(\mathrm{ad}\rho(h))^{n}\rho(x)\\
&=\mathrm{Exp}_{\mathfrak{D}}\circ \mathrm{ad}_{\mathfrak{D}}(\rho(h))(\rho(x))=\mathrm{Ad}_{\mathfrak{D}}(exp_{\mathfrak{D}}\rho(h))\rho(x)\\
&=\mathrm{Ad}_{\mathfrak{D}}(\pi(exph))\rho(x)=\pi(exph)\rho(x)\pi(exp(-h))\\
&=\pi(exph)\rho(x)\pi(exph)^{-1}
\end{split}
&
\end{align*}
Therefore for every 
$x \in \mathfrak{g_{{0}}}$,
this relation is established.
\\
Thus the 4-tuple $(\pi,\mathcal{H}_{\check{f},U},\mathcal{D}_{\check{f},U},\rho  )$ satisfies the properties stated in Definition \ref{f1}.
\section{	Extension of positive definite local $\mathbb{Z}_{2}^{n}$-superfunctions}\label{e11}
In the last section, associated to a local positive $\mathbb{Z}_{2}^{n}$-superfunction, a pre-representation, say $(\pi,\mathcal{H}_{\check{f},U},\mathcal{D}_{\check{f},U},\rho  )$ has been constructed. By Stability Theorem, one can extend this to a smooth unitary representation $(\pi,\rho^{\pi},\mathcal{H}_{\check{f},U})$ of $ \mathbb{Z}_2^n $-Harish-Chandra pair $ 	(G_0,\mathfrak{g}_{\mathbb{C}}) $. In this section, by using this representation, we construct an extension of the local positive definite superfunction $ f_{UU^{-1}} $.
\subsection{Constructing a global positive definite $\mathbb{Z}_{2}^{n}$-superfunction associated with the smooth unitary representation.}
For $1=(1,1)\in S$, we define matrix coefficients of the smooth unitary representation $(\pi,\rho^{\pi},\mathcal{H}_{\check{f},U})$ introduced in second step as follows. Set 
$v:=K_1$
\begin{align*}
&& \check{h}:\mathcal{S} \rightarrow \mathbb{C} \qquad \check{h}(g,D)=\phi_{v,v}(g,D):=( \pi(g) \rho^{\pi}(D) v,v)&&
\end{align*}
We show that map
$ h $
with the following property is a
positive definite smooth 	$\mathbb{Z}_{2}^{n}$-superfunctions of
$ (G_0,\mathfrak{g}_{\mathbb{C}}) $.
\begin{align*}
&\begin{split}
h: \mathfrak{U}(\mathfrak{g}_{\mathbb{C}}) \rightarrow C^{\infty}(G_0,\mathbb{C}),\qquad h(D)(g):=\check{h}(g,D)&
\end{split}
&
\end{align*}
To this end, it is sufficient to prove the statements (i) and (ii).
\begin{itemize}
	\item [(i)]
	$L_{x}(h(D))(g)=h(xD)(g)$
	for every 
	$x \in \mathfrak{g_{0}}$, 
	$D \in \mathfrak{U}(\mathfrak{g}_{\mathbb{C}})$
	and
	$g \in G_0$.\\
	In order to prove this equality, we use relations
	$ L_{x}(h(D))(g)=\lim\limits_{t \rightarrow 0}\frac{1}{t}\left(h(D)(ge^{tx})-h(D)(g)\right) $
	and
	$ L_{x}\rho^{\pi}(D)v=\rho^{\pi}(xD)v $.
	
	\begin{align*}
	\begin{split}
	L_{x}(h(D))(g)&=\lim\limits_{t \rightarrow 0}\frac{1}{t}\left(h(D)(ge^{tx})-h(D)(g)\right)\\&
	=\lim\limits_{t \rightarrow 0}\frac{1}{t}\left(\phi_{v,v}(ge^{tx},D)-\phi_{v,v}(g,D)\right)\\
	&=\lim\limits_{t \rightarrow 0}\frac{1}{t}\left((\pi(ge^{tx})\rho^{\pi}(D)v,v)-(\pi(g)\rho^{\pi}(D)v,v)\right)\\
	&=\lim\limits_{t \rightarrow 0}\frac{1}{t}\left(\pi(ge^{tx})\rho^{\pi}(D)v-\pi(g)\rho^{\pi}(D)v,v\right)\\
	&=\left(\lim\limits_{t \rightarrow 0}\frac{1}{t}\pi(ge^{tx})\rho^{\pi}(D)v-\pi(g)\rho^{\pi}(D)v,v\right)\\
	&=\left(L_{x}\pi(g)\rho^{\pi}(D)v,v\right)
	=\left(\pi(g)L_{x}\rho^{\pi}(D)v,v\right)\\
	&=\left(\pi(g)\rho^{\pi}(xD)v,v\right)=\phi_{v,v}(g,xD)=h(xD)(g).
	\end{split}
	&
	\end{align*}
	\item[(ii)]
	The mapping
	$ \phi_{n,v,v} $
	for every 
	$n \geq 0$
	is smooth where 
	$ \rho^{\pi}(x_1,...,x_n)v:=\rho^{\pi}(x_1)...\rho^{\pi}(x_n)v $.\\ Set
	\begin{align*}
	\begin{split}
	\phi_{n,v,v}:\mathfrak{g}^{n} \times G \rightarrow \mathbb{C},\quad \phi_{n,v,v}(x_1,...,x_n,g)&:=\left(\pi(g)\rho^{\pi}(x_1,...,x_n)v,v\right)\\&=\left(\rho(x_1,...,x_n)v,\pi(g^{-1})v\right)
	\end{split}
	&
	\end{align*}
	The smoothness of the first component is due to the continuity of the $n$-line mapping
	$(x_1,,...,x_n) \mapsto \rho(x_1,...,x_n)v$ and given that 	$v \in \mathcal{H}_{\check{f},U}$ then the second component meaning that the mapping
	$g \mapsto \pi(g^{-1})v$ is smooth.	
\end{itemize}
Thus mapping $\check{h}$ is a matrix coefficient of the above unitary representation of $ (G_0,\mathfrak{g}_{\mathbb{C}}) $, hence it relates to a positive definite function in $ (G_0,\mathfrak{g}_{\mathbb{C}}) $.
\subsection{ Extension 	of positive definite local $\mathbb{Z}_{2}^{n}$-superfunction}
Here we state our main theorem as follows:
\begin{theorem}\label{e8}
	Let $ G_0 $ be a connected Lie group and $ V $ be an open neighborhood of identity element, then each local positive definite $\mathbb{Z}_{2}^{n}$-superfunction say $f_{V} \in \mathrm{Hom}_{\mathfrak{g}_{0}}(\mathfrak{U}(\mathfrak{g}_{\mathbb{C}}),C_{G_{0}}^{\infty}(V))$, if belongs to category $C_{ub}^{\infty}$, has a positive definite extension to all of the $\mathbb{Z}_{2}^{n}$-Lie supergroup, i.e. there exists a positive definite $\mathbb{Z}_{2}^{n}$-superfunction $ \tilde{f}\in Hom(\mathfrak{U}(\mathfrak{g}_{\mathbb{C}}), C^{\infty}(G_0)) $ and a open neighborhood of identity element $ U \subset V $ such that $ \tilde{f}|_U=f_V|_U $.
\end{theorem}
\proof
	Here we show that the global positive definite $\mathbb{Z}_{2}^{n}$-superfunction $ h $, introduced at the start point of  the last subsection is an extension of $f_{UU^{-1}}$. For this, we show that 
	$\check{h}\arrowvert_{ \mathcal{S}_U}=\check{f}$, i.e. for
	$(p,D) \in S_{U}$, one has
	$\check{h}(p,D)=\check{f}(p,D)$.\\
	Let 
	$p \in U \subset G_0$ then there exist tangent vectors
	$h_1,h_2,...,h_k \in \mathfrak{g_{{0}}}$
	such that 
	$p=exph_1...exph_k$.
	Define
	\begin{align*}
	&&\pi(p):=\pi(exph_1)...\pi(exph_k)&&
	\end{align*}
	If
	$p=exph_1$.
	(Here we display $f_{UU^{-1}}$ with $f$ for convenience)
	\begin{align*}
	\begin{split}
	\check{h}(p,D)&=(\pi(p)\rho(D)K_1,K_1)=(\pi(p)K_{(1,D^*)},K_{1})=(\pi(exph_{1})K_{(1,D^*)},K_1)\\&=(\tilde{\rho}(h_1)K_{(1,D^*)},K_1)
	=(e^{\overline{\rho(h_1)}}K_{(1,D^*)},K_1)=(\sum\limits_{r=0}^{\infty}\frac{\rho(h_1)^{r}}{r!}K_{(1,D^*)},K_{1})\\&=\sum\limits_{r=0}^{\infty}\dfrac{1}{r!}(\rho(h_1)^rK_{(1,D^*)},K_{1})
	=\sum\limits_{r=0}^{\infty}\dfrac{1}{r!}(K_{(1,D^*)(1,{h_{1}^{*}}^{r})},K_{1})\\&=\sum\limits_{r=0}^{\infty}\dfrac{1}{r!}(K_{(1,D^*{h_{1}^{*}}^{r})},K_{1})=\sum\limits_{r=0}^{\infty}\dfrac{1}{r!}\check{f}(1,h_{1}^{r}D)
	=\sum\limits_{r=0}^{\infty}\dfrac{1}{r!}
	f(h_{1}^{r}D)(1)
	\end{split}
	&
	\end{align*}
	since superfunction 
	$f$ belongs to
	$C_{ub}^{\infty}(V)$, then series 
	$\sum\limits_{r=0}^{\infty}\dfrac{1}{r!}
	f(h_{1}^{r}D)(1)$ is convergent.
	$$\sum\limits_{r=0}^{\infty}\dfrac{1}{r!}
	f(h_{1}^{r}D)(1)=\sum\limits_{r=0}^{\infty}\dfrac{1}{r!}D_{h_{1}}^{r}f(D)(1)=\sum\limits_{r=0}^{\infty}\dfrac{1}{r!}\frac{d^r}{dt^r}f(D)(expth_{1})|_{t=0}$$
	Since the following equation holds
	$$\sum\limits_{r=0}^{\infty}\dfrac{1}{r!}\frac{d^r}{dt^r}f(D)(expth_{1})|_{t=0}t^r=f(D)(expth_1)$$
	By setting
	$t=1$, we have 
	$$ \sum\limits_{r=0}^{\infty}\dfrac{1}{r!}\frac{d^r}{dt^r}f(D)(expth_{1})|_{t=0 }=f(D)(exph_{1})$$
	Then
	$$\check{h}(p,D)=f(D)(exph_{1})=\check{f}(exph_{1},D)=\check{f}(p,D).$$
	Therefore 
	$\check{f}$
	and 
	$\check{h}$ are equal
	on 
	$\mathcal{S}_U$. By induction, we get the result.

{\footnotesize\begin{thebibliography}{22}
		\bibitem{Akutowicz1959}
		Akutovich, E. J., On Extrapolating a Positive Definite Function from a
		Finite Interval, Math. Scand. 7 (1959), 157-169.
		\bibitem{Berg2008}
		Berg, Christian. "Stieltjes-Pick-Bernstein-Schoenberg and their connection to complete monotonicity." Positive Definite Functions: From Schoenberg to Space-Time Challenges (2008), 15-45.
		\bibitem{Bruzual2001}
		Bruzual, Ramon, and Marisela Dominguez. "On the local operator commutation relations and extensions of locally defined positive definite functions on the Heisenberg group." Journal of mathematical analysis and applications 253.1 (2001), 215-223.
		\bibitem{carmeli2011mathematical}
		Carmeli, Claudio, Lauren Caston, and Rita Fioresi. Mathematical foundations of supersymmetry. Vol. 15. European Mathematical Society, 2011.
		\bibitem{Conway1985}
		Conway, John B. A Course in Functional Analysis. Graduate Texts in
		Mathematics 96. Springer, 2nd edition, 1985.
		\bibitem{CTGJPN}
		Covolo, Tiffany,  Grabowski, Janusz, and  Poncin, Norbert. "$\mathbb {Z} _2^ n $-Supergeometry I: Manifolds and Morphisms." arXiv preprint arXiv:1408.2755, 2014.
		\bibitem{covolo2016differential}
		Covolo, Tiffany, Stephen Kwok, and Norbert Poncin. "Differential calculus on $\mathbb {Z}^ n_2 $-supermanifolds." arXiv preprint arXiv:1608.00949 (2016).
		\bibitem{Devinatz1959}
		Devinatz, A. On the Extensions of Positive Definite Functions, Acta
		Math., 102 (1959), 109-134.
		\bibitem{Friedrich1991}
		Friedrich, Jurgen. Extension of positive definite functions on lie groups. In
		Seminar Sophus Lie, Citeseer., volume 1  (1991), 225-233.
		\bibitem{HB2015}
		Hall, Brian C. Lie Groups, Lie Algebras, and Representations: An Elementary Introduction, Graduate Texts in Mathematics, vol. 222 (2nd ed.), Springer, 2015.
		\bibitem{JorgensenPedersen2016}
		Jorgensen, Palle, Steen Pedersen, and Feng Tian. Extensions of Positive Definite Functions: Applications and Their Harmonic Analysis. Vol. 2160. Springer, 2016.
		\bibitem{Jorgensen1989}
		Jorgensen, Palle ET. "Positive definite functions on the Heisenberg group." Mathematische Zeitschrift 201.4 (1989), 455-476.
		\bibitem{Jorgensen1990}
		|,"Extensions of positive definite integral kernels on the Heisenberg group." Journal of functional analysis 92.2 (1990), 474-508.
		\bibitem{Jorgensen1991}
		|,"Integral representations for locally defined positive definite functions on Lie groups." International Journal of Mathematics 2.03 (1991), 257-286.
		\bibitem{Kaniuth2004}
		Kaniuth, Eberhard. "Extension of positive definite functions from subgroups of nilpotent locally compact groups." Proceedings of the American Mathematical Society 132.3 (2004): 865-874.
		\bibitem{Krein1940}
		Krein, Mark G. "Sur le probleme du prolongement des fonctions hermitiennes positives et continues." CR (Doklady) Acad. Sci. URSS (NS) 26.1 (1940), 17-22.
		\bibitem{Livshic1945}
		Livshic, M. S., On Self-adjoint and Skew-adjoint Extensions of Symmetric
		Operators, Moscow, 1945, (Russian).
		\bibitem{Manin1988}	
		Manin, Yuri I. Gauge Field Theory and Complex Geometry, Grundlehren der Mathematischen
		Wissenschaften [Fundamental Principles of Mathematical Sciences], vol. 289, Springer-Verlag,
		Berlin, 1988, Translated from the Russian by N. Koblitz and J. R. King. MR MR954833
		(89d:32001).
		\bibitem{MSNKHSH}
		Merigon, Stephane, Neeb, Karl-Hermann, and Salmasian, Hadi. Categories
		of unitary representations of banach Lie supergroups and restriction functors.
		Pacific Journal of Mathematics, 257(2):431-469, 2012.
		\bibitem{MMSH}
		Mohammadi, Mohammad and Salmasian, Hadi. The Gelfand-Naimark-
		Segal construction for unitary representatins of $\mathbb{Z}_{2}^{n}$ -graded lie supergroups.
		Banach Center Publications, 113 (2017), 263-274.
		\bibitem{mohammadi2021construction}
		Mohammadi, Mohammad, and Saad Varsaie. "On the construction of supergrassmannians as homogeneous superspaces." Asian-European Journal of Mathematics 14.04 (2021): 2150053.
		\bibitem{Nee11}
		Neeb, Karl-H. "On analytic vectors for unitary representations of infinite dimensional Lie groups." In Annales de l'institut Fourier, vol. 61, no. 5, pp. 1839-1874. 2011.
		\bibitem{NKHSH}
		Neeb, Karl-Hermann and Salmasian, Hadi. Differentiable vectors and unitary
		representations of Frechet-Lie supergroups. Math. Z., 275(1-2):419-451, 2013.
		\bibitem{NKH}
		Neeb, Karl-Hermann. On differentiable vectors for representations of infinite dimensional Lie groups. Journal of Functional Analysis 259.11: 2814-2855,2010.
		\bibitem{NeebSalmasian2013}
		Neeb, Karl-Hermann, and Hadi Salmasian. "Positive definite superfunctions and unitary representations of Lie supergroups." Transformation Groups 3.18 (2013), 803-844.
		\bibitem{NK}
		Nishiyama, Kyo. "Oscillator representations for orthosymplectic algebras." Journal of Algebra 129.1: 231-262, 1990.
		\bibitem{Raikov1940}
		Raikov, DA. Sur les fonctions positivement definies. Dokl. Akad. Nauk
		SSSR, 26 (1940), 860-865.
		\bibitem{RMSB}
		Reed, Michael and Simon Barry. Methods of modern mathematical physics. Functional
		analysis, volume 1. Academic Press, rev and enl edition, 1981.
		\bibitem{Rudin1963}
		Rudin, Walter. "The extension problem for positive-definite functions." Illinois J. Math. 7.1 (1963): 532-539.
		\bibitem{Rudin1970}
		Rudin, Walter. "An extension theorem for positive-definite functions." Duke Mathematical Journal 37.1 (1970): 49-53.
		\bibitem{Sasvari1944}
		Sasvari, Z. Positive definite and definitizable functions, akademie verlag,
		berlin, 1994.
		\bibitem{Stewart1976}
		Stewart, James. "Positive definite functions and generalizations, an historical survey." The Rocky Mountain Journal of Mathematics 6.3 (1976), 409-434.
\end{thebibliography}}
\end{document}